\documentclass[10pt]{article}

\usepackage{amssymb,amsthm,amsmath,amsfonts}
\usepackage{graphicx}
\usepackage{booktabs,url}
\usepackage{algorithm,algorithmic}
\usepackage[usenames]{color}
\graphicspath{{./pngfig/}}

\title{A Variational Bayesian Method to Inverse Problems with Impulsive Noise}
\date{\today}
\author{
Bangti Jin\footnote{Department of Mathematics and Institute for Applied
Mathematics and Computational Sciences, Texas A\&M University,
College Station 77843-3368, TX, USA. (btjin@math.tamu.edu).}}
\begin{document}
\maketitle

\begin{abstract}
We propose a novel numerical method for solving inverse problems subject to impulsive
noises which possibly contain a large number of outliers. The approach is of Bayesian
type, and it exploits a heavy-tailed $t$ distribution for data noise to achieve
robustness with respect to outliers. A hierarchical model with all hyper-parameters
automatically determined from the given data is described. An algorithm of variational
type by minimizing the Kullback-Leibler divergence between the true posteriori
distribution and a separable approximation is developed. The numerical method is
illustrated on several one- and two-dimensional linear and nonlinear inverse problems
arising from heat conduction, including estimating boundary temperature, heat flux and
heat transfer coefficient. The results show its robustness to outliers and the fast and
steady convergence of the algorithm.
\\
\textbf{key words}: impulsive noise, robust Bayesian, variational method, inverse
problems
\end{abstract}

\section{Introduction}

We are interested in Bayesian approaches for inverse problems subject to impulsive
noises. Bayesian inference provides a principled framework for solving diverse inverse
problems, and has demonstrated distinct features over deterministic techniques, e.g.,
Tikhonov regularization. Firstly, it can yield an ensemble of plausible solutions
consistent with the given data. This enables quantifying the uncertainty of a specific
solution, e.g., with credible intervals. In contrast, deterministic techniques generally
content with singling out one solution from the ensemble. Secondly, it provides a
flexible regularization since hierarchical modeling can partially resolve the nontrivial
issue of choosing an appropriate regularization parameter. It is known that the
underlying mechanism is balancing principle \cite{JinZou:2009ip}. Thirdly, it allows
seamlessly integrating structural/multiscale features of the problem through careful
prior modeling. Therefore, it has attracted attention in a wide variety of applied
disciplines, e.g., geophysics \cite{Tarantola:2005,SambridgeMosegaard:2002}, medical
imaging
\cite{HebertLeahy:1989,KaipioKolehmainenSomersaloVauhkonen:2000,AchimBezerianosTsakalides:2001}
and heat conduction
\cite{Emery:2002,WangZabaras:2004ijhmt,WangZabaras:2005ip,Emery:2009}, see also
\cite{MarzoukNajmRahn:2007,Koutsourelakis:2009,MaZabaras:2009,MarzoukXiu:2009} for other
applications. For an overview of methodological developments, we refer to the monographs
\cite{Tarantola:2005,KaipioSomersalo:2005}.

Amongst existing studies on Bayesian inference for inverse problems, the Gaussian noise
model has played a predominant role. This is often justified by appealing to central
limit theorem. The theorem asserts that the normal distribution is a suitable model for
data that are formed as the sum of a large number of independent components. Even in the
absence of such justifications, this model is still preferred due to its
computational/analytical conveniences, i.e., it allows direct computation of the
posterior mean and variance and easy exploration of the posterior state space (for linear
models with Gaussian priors). A well acknowledged limitation of the Gaussian model is its
lack of robustness against the outliers, i.e., data points that lie far away from the
bulk of the data, in the observations: A single aberrant data point can significantly
influence all the parameters in the model, even for these with little substantive
connection to the outlying observations \cite[pp. 443]{GelmanCarlinSternRubin:2004}.

However, it is clear that not all real-world data can be adequately described by the
Gaussian model. For example, Laplacian noises can arise when acquiring certain signals
\cite{AllineyRuzinsky:1994}, and salt-and-pepper noises are very common in natural
images/signals due to faulty memory location and transmission in noisy channels
\cite{ChanHoNikolova:2005}. The impulsive nature of these noises is reflected by the
heavy tail of their distributions and thus the presence of, possibly of a significant
amount, outliers in the data. Physically, such noises arise from uncertainties in
instrument calibration, physical limitations of acquisition devices and experimental
(operation) conditions. Due to its lack of robustness, an inadvertent adoption of the
Gaussian model can seriously compromise the accuracy of the estimate, and consequently
does not allow full extraction of the information provided by the data.

This calls for methods that are robust to the presence of outliers. There are several
ways to derive robust estimates. One classical approach is to first identify the outliers
with noise detectors, e.g., by adaptive media filter, and then to perform
inversion/reconstruction on the data set with outliers excluded
\cite{GelmanCarlinSternRubin:2004,ChanHoNikolova:2005}. The success of such procedures
relies crucially on the reliability of the noise detector. However, it can be highly
nontrivial to accurately identify all outliers, especially in high dimensions, and
mis-identification can adversely affect the quality of subsequent inversion. This
necessitates developing systematic strategies for handling impulsive noises, which can be
achieved by modeling the outliers explicitly with a heavy-tailed noise distribution. The
Student's $t$ and Laplace distributions are two most popular choices. The use of the
$t$-distribution in robust Bayesian analysis is well recognized, see
\cite{LangeLittleTaylor:1989} for its usage in statistical contexts.  In
\cite{DempsterLairdRubin:1977}, the application of the EM algorithm to $t$ models was
shown. Recently, Tipping and Lawrence \cite{Tipping:2005} developed a robust Bayesian
interpolation with the $t$-distribution. Alternatively, the Laplace distribution may be
employed, see, e.g., \cite{Gao:2008} for robust probabilistic principal component
analysis.

This paper studies the potentials of one robust Bayesian formulation for inverse problems
subject to impulsive noises. Impulsive noise has received some recent attention in
deterministic inversion, see \cite{ClasonJinKunisch:2010sii,Clason:2009b} and references
therein, but not in the Bayesian framework. The salient features of the proposed approach
include uncertainty quantification of the computed solution, robustness to data outliers,
and general applicability to both linear and nonlinear inverse problems. Therefore, it
complements the developments of robust formulations in the framework of deterministic
inverse problems \cite{ClasonJinKunisch:2010sii,Clason:2009b}. As to the numerical
exploration of the Bayesian model, we capitalize on the variational method developed in
machine learning \cite{JordanGhahramaniJaakkolaSaul:1999,Attias:2000,Beal:2003} for
approximate inference, and thus achieve reasonable computational efficiency. The
application of variational Bayesian formulations to inverse problems, especially
nonlinear ones, is of relatively recent origin \cite{SatoTaku:2004,JinZou:2010jcp}, and
their robust counterparts seem largely unexplored.

The rest of the paper is structured as follows. A hierarchical formulation based on the
$t$ distribution for the noise is derived in Section \ref{sect:bayes}. The variational
method for numerically exploring the posterior is described in Section \ref{sect:vb}, and
two algorithms are developed for linear and nonlinear inverse problems, respectively. In
Section \ref{sec:exp}, numerical results for four benchmark inverse problems arising in
heat transfer are presented to illustrate the features of the formulation and the
convergence behavior of the algorithms. 

\section{Hierarchical Bayesian inference}\label{sect:bayes}
In this section, we formulate the hierarchical Bayesian model for inverse problems
subject to impulsive noises. The focus is on the noise model and hyper-parameter
treatment.

We consider for the following finite-dimensional linear inverse problem
\begin{equation}\label{eqn:inv}
\mathbf{K(u)}=\mathbf{y},
\end{equation}
where $\mathbf{K}:\mathbb{R}^m\to\mathbb{R}^n$, $\mathbf{u}\in\mathbb{R}^m$ and
$\mathbf{y}\in\mathbb{R}^n$ represent the (possibly nonlinear) forward model, the
sought-for solution and noisy observational data, respectively.

In Bayesian formalism, the likelihood function $p(\mathbf{y}|\mathbf{u})$ incorporates
the information contained in the data $\mathbf{y}$, and it is dictated by the noise
model. Let the given data $\mathbf{y}$ be subjected to additive noises, i.e.,
\begin{equation*}
\mathbf{y}=\mathbf{y}^\dagger+\boldsymbol{\zeta},
\end{equation*}
where $\boldsymbol{\zeta}\in \mathbb{R}^n$ is a random vector corrupting the exact data
$\mathbf{y}^\dagger$. In practice, a Gaussian distribution on each component $\zeta_i$ is
customarily assumed. The validity of this assumption relies crucially on being not
heavy-tailed and the symmetry of the distribution, and the violation of either condition
may render the resulting Bayesian model invalid and inappropriate, which may seriously
compromise the accuracy of the posteriori estimate.

In practice, data outliers can arise due to, e.g., erroneous recording and transmission
in noisy channels, which makes the Gaussian model unsuitable. Following the interesting
works \cite{LangeLittleTaylor:1989,Geweke:1993,Tipping:2005}, we choose to model the
outliers explicitly by a heavy-tailed distribution, and this yields a seamless and
systematic framework for treating impulsive noises. We focus on the $t$ model
\cite{GelmanCarlinSternRubin:2004}, where the noises $\zeta_i$ are independent and
identically distributed according to a centered $t$ distribution, i.e.,
\begin{equation*}
  p(\zeta_i;\nu,\sigma)=\frac{\Gamma(\frac{\nu+1}{2})}{\Gamma(\frac{\nu}{2})\sqrt{\pi\nu}\sigma}
  \left\{1+\frac{1}{\nu}\left(\frac{\zeta_i}{\sigma}\right)^2\right\}^{-\frac{\nu+1}{2}},
\end{equation*}
where $\nu$ is a degree of freedom parameter, $\sigma$ is a scale parameter
\cite{GelmanCarlinSternRubin:2004}, and $\Gamma(\cdot)$ is the standard Gamma function.
Consequently, the likelihood function $p(\mathbf{y}|\mathbf{u})$ is given by
\begin{equation}\label{eqn:like}
p(\mathbf{y}|\mathbf{u})=\left(\frac{\Gamma(\frac{\nu+1}{2})}{\Gamma(\frac{\nu}{2})\sqrt{\pi\nu}\sigma}\right)^n
  \prod_{i=1}^n\left\{1+\frac{1}{\nu}\left(\frac{|(\mathbf{K(u)-y})_i|}{\sigma}\right)^2\right\}^{-\frac{\nu+1}{2}},
\end{equation}
where the subscript $i$ denotes the $i$th entry of a vector.

In Bayesian formalism, structural prior knowledge about the unknown $\mathbf{u}$ is
encoded in the prior distribution $p(\mathbf{u})$. Here we focus on the following simple
random field
\begin{equation}\label{eqn:mrf}
p(\mathbf{u}|\lambda)=C
\lambda^{\frac{s}{2}}\exp\left(-\frac{\lambda}{2}\|\mathbf{Lu}\|_2^2\right),
\end{equation}
where $\|\cdot\|_2$ denotes the Euclidean norm and $C$ is a normalizing constant. The
matrix $\mathbf{L}\in\mathbb{R}^{s\times m}$, whose rank is $s$, encodes the structural
interactions between neighboring sites, and $\lambda$ is a scaling parameter.

The hyper-parameters $\nu$, $\sigma$ and $\lambda$ in the likelihood
$p(\mathbf{y}|\mathbf{u})$ and the prior $p(\mathbf{u}|\lambda)$ play the crucial role of
regularization parameters in classical regularization \cite{JinZou:2009ip}. Bayesian
formalism resolves the issue through hierarchical modeling and determines them
automatically from the data $\mathbf{y}$.  A standard practice to select priors for
hyper-parameters is to use conjugate priors. For the parameter $\lambda$, the conjugate
prior is a Gamma distribution $G(t;\alpha,\beta)$, which is defined by
\begin{equation}\label{eqn:gamdist}
G(t;\alpha,\beta)=\frac{\beta^\alpha}{\Gamma(\alpha)}
t^{\alpha-1}e^{-\beta t},
\end{equation}
where $\alpha$ and $\beta$ are nonnegative constants. The parameters $\nu$ and $\sigma$
do not admit easy conjugate form, and one may opt for the maximum likelihood approach
when appropriate.

According to Bayes' rule, the posterior $p(\mathbf{u},\lambda|\mathbf{y})$ is related to
the data $\mathbf{y}$ by
\begin{equation*}
p(\mathbf{u},\lambda|\mathbf{y})=\frac{p(\mathbf{y}|\mathbf{u})p(\mathbf{u}|\lambda)p(\lambda)}{\int\int
p(\mathbf{y}|\mathbf{u})p(\mathbf{u}|\lambda)p(\lambda)d\mathbf{u}d\lambda}.
\end{equation*}
Upon ignoring the (unimportant) normalizing constant $p(\mathbf{y})=\int\int p(\mathbf{y}
|\mathbf{u})p(\mathbf{u}|\lambda)p(\lambda)d\mathbf{u} d\lambda$, the posterior
$p(\mathbf{u},\lambda|\mathbf{y})$ may be simply evaluated as
\begin{eqnarray}\label{eqn:ppdf}
p(\mathbf{u},\lambda|\textbf{y})\propto \left(\frac{\Gamma(\frac{\nu+1}{2})}{\Gamma(\frac{\nu}{2})\sqrt{\pi\nu}\sigma}\right)^n
  \prod_{i=1}^n\left\{1+\frac{1}{\nu}\left(\frac{|(\mathbf{K(u)-y})_i|}{\sigma}\right)^2\right\}^{-\frac{\nu+1}{2}}
  \cdot\lambda^{\frac{s}{2}}e^{-\frac{\lambda}{2}\|\mathbf{Lu}\|_2^2}\cdot \lambda^{\alpha_0-1}e^{-\beta_0 \lambda},
\end{eqnarray}
where $(\alpha_0,\beta_0)$ is the parameter pair of the Gamma distribution for $\lambda$.

The posterior state space $p(\mathbf{u},\lambda|\mathbf{y})$ is often high dimensional,
and thus it can only be numerically explored. In Section \ref{sect:vb}, we shall develop
an efficient variational method for its approximate inference.

\section{Variational approximation}\label{sect:vb}
In this section, we describe a variational method for efficiently constructing an
approximation to the posterior distribution \eqref{eqn:ppdf}. It can deliver point
estimates together with uncertainties for both the solution $\mathbf{u}$ and the
hyper-parameter $\lambda$. There are three major obstacles in getting a faithful
approximation:
\begin{itemize}
 \item[(a)] nongaussian likelihood ($t$ instead of Gaussian distribution),
 \item[(b)] statistical dependency between $\mathbf{u}$ and $\lambda$, and
 \item[(c)] possible nonlinearity of the forward mapping $\mathbf{K}$.
\end{itemize}
To circumvent these obstacles, we shall make use of three ideas: scale-mixture
representation of the $t$ distribution, variational (separable) approximation for
decoupling dependency, and recursive linearization for resolving nonlinearity.

First we describe the scale-mixture representation. In the posterior \eqref{eqn:ppdf},
the $t$ likelihood makes it hard to find or define a good approximation. Fortunately, it
can be represented as follows \cite[pp. 446]{GelmanCarlinSternRubin:2004}
\begin{equation*}
  \begin{aligned}
    p(\zeta_i|\nu,\sigma)&=\frac{\Gamma(\frac{\nu+1}{2})}{\Gamma(\frac{\nu}{2})\sqrt{\pi\nu}\sigma}
    \left\{1+\frac{1}{\nu}\left(\frac{\zeta_i}{\sigma}\right)^2\right\}^{-\frac{\nu+1}{2}}\\
    &=\int_0^\infty\sqrt{\frac{w_i}{2\pi}}e^{-\frac{w_i}{2}\zeta_i^2}p(w_i;\nu,\sigma)dw_i,
  \end{aligned}
\end{equation*}
where the density $p(w_i)$ is given by $ p(w_i;\nu,\sigma)=\frac{(\frac{\nu\sigma^2}{2}
)^{\frac{\nu}{2}}}{\Gamma(\frac{\nu}{2})}w^{\frac{\nu}{2}-1}
e^{-\frac{\nu\sigma^2}{2}w}=G\left(w_i;\tfrac{\nu}{2},\tfrac{\nu\sigma^2}{2}\right)$,
c.f. \eqref{eqn:gamdist}. To simplify the expression, we introduce two independent
variables $\alpha_1$ and $\beta_1$ by
\begin{equation*}
\alpha_1=\frac{\nu}{2}\quad \mbox{and}\quad \beta_1=\frac{\nu\sigma^2}{2},
\end{equation*}
and work with the parameters $\alpha_1$ and $\beta_1$ hereon. Then we have the following
succinct formula
\begin{equation*}
  p(\zeta_i|\alpha_1,\beta_1)=\int_0^\infty\sqrt{\frac{w_i}{2\pi}}e^{-\frac{w_i}{2}\zeta_i^2}p(w_i;\alpha_1,\beta_1)dw_i,
\end{equation*}
with $p(w_i;\alpha_1,\beta_1)=G\left(w_i;\alpha_1,\beta_1\right)$. Therefore, the $t$
model is a mixture (average) of an infinite number of Gaussians of varying precisions
$w_i$, with the mixture weight $w_i$ specified by the Gamma distribution
$p(w_i;\alpha_1,\beta_1)$. The representation also explains its heavy tail: for small
$w_i$, the random variable $\zeta_i$ follow a Gaussian distribution with a large
variance, and thus the realizations are likely to take large values, which behaves more
or less like outliers. By means of scale mixture, we have introduced an extra variable,
but effectively converted a $t$ distribution into a Gaussian distribution. In sum, we
have arrived at the following augmented posterior
\begin{equation}\label{eqn:appdf}
p(\mathbf{u},\mathbf{w},\lambda|\mathbf{y})\propto
|\mathbf{W}|^\frac{1}{2}e^{-\frac{1}{2}\|\mathbf{K(u)}-\mathbf{y}\|_\mathbf{W}^2}
\cdot p(\mathbf{w};\alpha_1,\beta_1)\cdot\lambda^{\frac{s}{2}}e^{-\frac{\lambda}{2}\|\mathbf{Lu}\|_2^2}
\cdot \lambda^{\alpha_0-1}e^{-\beta_0\lambda},
\end{equation}
where $\mathbf{w}\in\mathbb{R}^n$ is an auxiliary random vector following the Gamma
distribution, i.e.,
\begin{equation*}
  p(\mathbf{w};\alpha_1,\beta_1) = \prod_{i=1}^n G(w_i;\alpha_1,\beta_1)=G(\mathbf{w};\alpha_1,\beta_1),
\end{equation*}
$\mathbf{W}$ is a diagonal matrix with diagonal $\mathbf{w}$, and the weighted norm
$\|\cdot\|_\mathbf{W}$ is defined by
$\|\mathbf{v}\|_\mathbf{W}^2=\mathbf{v}^\mathrm{T}\mathbf{Wv}$. The posterior
$p(\mathbf{u},\mathbf{w},\lambda|\mathbf{y})$ is computationally more amenable with the
variational method.

Next we describe the variational method for approximately exploring the posterior
\eqref{eqn:appdf} in case of a linear operator $\mathbf{K}$, i.e.,
$\mathbf{K}(\mathbf{u})=\mathbf{Ku}$. The derivations here follow closely
\cite{JinZou:2010jcp,Tipping:2005}. An approximation can be derived as follows. One first
transforms it into an equivalent optimization problem using the Kullback-Leibler
divergence and then obtains an approximation by solving the optimization problem
inexactly. The divergence $D_{KL}(q(\mathbf{u},\mathbf{w},\lambda)|
p(\mathbf{u},\mathbf{w},\lambda|\mathbf{y}))$ between two densities
$q(\mathbf{u},\mathbf{w},\lambda)$ and $p(\mathbf{u},\mathbf{w},\lambda|\mathbf{y})$ is
defined by
\begin{equation*}
\begin{aligned}
D_{KL}(q(\mathbf{u},\mathbf{w},\lambda)|p(\mathbf{u},\mathbf{w},
\lambda|\mathbf{y}))&=\int\int\int q(\mathbf{u},\mathbf{w},\lambda)
\ln\frac{q(\mathbf{u},\mathbf{w},\lambda)
}{p(\mathbf{u},\mathbf{w},\lambda,\mathbf{y})} d\mathbf{u}d\mathbf{w}d\lambda+\log p(\mathbf{y}),
\end{aligned}
\end{equation*}
where $p(\mathbf{y})$ is a normalizing constant. Since the divergence $D_{KL}$ is
nonnegative and vanishes if and only if $q$ coincides with $p$, minimizing $D_{KL}$
effectively transforms the problem into an equivalent optimization problem. We shall
minimize the following functional, which is also denoted by $D_{KL}$
\begin{equation}\label{eqn:kld}
D_{KL}(q(\mathbf{u},\mathbf{w},\lambda)|p(\mathbf{u},\mathbf{w},\lambda|\mathbf{y}))=\int\int\int
q(\mathbf{u},\mathbf{w},\lambda) \ln\frac{q(\mathbf{u},\mathbf{w},\lambda)
}{p(\mathbf{u},\mathbf{w},\lambda,\mathbf{y})} d\mathbf{u}d\mathbf{w}d\lambda.
\end{equation}

However, directly minimizing $D_{KL}$ is still intractable since the posterior
$p(\mathbf{u}, \mathbf{w}, \lambda|\mathbf{y})$ is not available in closed form. We
impose a separability (conditionally independence) condition for the posterior
distributions of $\mathbf{u}$, $\mathbf{w}$ and $\lambda$ to arrive at a tractable
approximation, i.e.,
\begin{equation}\label{eqn:fac}
q(\mathbf{u},\mathbf{w},\lambda)=q(\mathbf{u})q(\mathbf{w})q(\lambda).
\end{equation}

\begin{algorithm}
\caption{Variational approximation for linear models $\mathbf{K}$}\label{alg:va}
\begin{algorithmic}[1]
 \STATE Set initial guess $q^0(\mathbf{w})$ and $q^0(\lambda)$;
 \FOR {$k=1,\dots,K$}
      \STATE Find $q^{k}(\mathbf{u})$ by
             \begin{equation*}
                  q^k(\mathbf{u})=\arg\min_{q(\mathbf{u})}
                  D_{KL}(q(\mathbf{u})q^{k-1}(\mathbf{w})q^{k-1}(\lambda)|p(\mathbf{u},\mathbf{w},\lambda|\mathbf{y}));
             \end{equation*}
      \STATE Find $q^{k}(\mathbf{w})$ by
             \begin{equation*}
                  q^k(\mathbf{w})=\arg\min_{q(\mathbf{w})}
                  D_{KL}(q^k(\mathbf{u})q(\mathbf{w})q^{k-1}(\lambda)|p(\mathbf{u},\mathbf{w},\lambda|\mathbf{y}));
             \end{equation*}
      \STATE Find $q^k(\lambda)$ by
             \begin{equation*}
                  q^k(\lambda)=\arg\min_{q(\lambda)}
                  D_{KL}(q^k(\mathbf{u})q^k(\mathbf{w})q(\lambda)|p(\mathbf{u},\mathbf{w},\lambda|\mathbf{y}));
             \end{equation*}
      \STATE Check a stopping criterion;
 \ENDFOR
\STATE Return approximation $q^k(\mathbf{u})q^k(\mathbf{w})q^k(\lambda)$.
\end{algorithmic}
\end{algorithm}

Under condition \eqref{eqn:fac}, an approximation can be computed by Algorithm
\ref{alg:va}. Each step of the algorithm can be further developed as follows. Setting the
first variation of $D_{KL}$ with respect to $q(\mathbf{u})$ to zero gives
\begin{equation*}
q^k(\mathbf{u})\propto\exp\left(E_{q^{k-1}(\mathbf{w})q^{k-1}(\lambda)}[\ln
p(\mathbf{u},\mathbf{w},\lambda,\mathbf{y})]\right).
\end{equation*}
It follows that $q^k(\mathbf{u})$ follows a Gaussian distribution with covariance
$\mathrm{cov}_{q^k(\mathbf{u})}$ and mean $\mathbf{u}_k$ given by
\begin{equation*}
\mbox{cov}_{q^k(\mathbf{u})}[\mathbf{u}]=\left[
\mathbf{K}^\mathrm{T}\mathbf{W}_k\mathbf{K}+\lambda_k
\mathbf{L}^\mathrm{T}\mathbf{L}\right]^{-1} \quad\mbox{and}\quad
\mathbf{u}_{k}:=E_{q^k(\mathbf{u})}[\mathbf{u}]=\mbox{cov}_{q^k(\mathbf{u})}
[\mathbf{u}]\mathbf{K}^\mathrm{T}\mathbf{W}_k\mathbf{y},
\end{equation*}
respectively, where $\lambda_k=E_{q^{k-1}(\lambda)}[\lambda]$ and
$\mathbf{W}_k=E_{q^{k-1}(\mathbf{w})}[\mathbf{W}]$, i.e.,
\begin{equation*}
q^k(\mathbf{u})=N(\mathbf{u};\mathbf{u}_k,[\mathbf{K}^\mathrm{T}\mathbf{W}_k\mathbf{K}+\lambda_k
\mathbf{L}^\mathrm{T}\mathbf{L}]^{-1}),
\end{equation*}
where $N$ refers to a normal distribution. Analogously, we can show that
$q^{k}(\mathbf{w})$ and $q^k(\lambda)$ take a factorized form, i.e.,
\begin{equation*}
\begin{aligned}
  q^{k}(\mathbf{w})&=G\left(\mathbf{w};\alpha_1+\tfrac{1}{2},\beta_1+\tfrac{1}{2}E_{q^k(\mathbf{u})}[|\mathbf{Ku}-\mathbf{y}|^2]\right),\\
  q^{k}(\lambda)&=G\left(\lambda;\alpha_0+\tfrac{s}{2},\beta_0+\tfrac{1}{2}E_{q^k(\mathbf{u})}
  [\|\mathbf{Lu}\|_2^2]\right).
\end{aligned}
\end{equation*}
Thus Steps 4 and 5 involve simply updating their respective parameter pairs.

There are several viable choices for the stopping criterion at Step 6. In practice, the
following two heuristics work well. One is to monitor the relative change of the inverse
solution $\mathbf{u}_k$. If the change between two consecutive iterations falls below a
given tolerance $tol$, i.e., $\|\mathbf{u}_{k}
-\mathbf{u}_{k-1}\|_2/\|\mathbf{u}_k\|_2\leq tol$, then the algorithm may stop. Another
is to monitor the variable $\lambda_k$. Numerically, we observe that the algorithm
converges reasonably fast and steadily.

We briefly remark on the choice of the pair $(\alpha_1,\beta_1)$. For our experiments in
Section \ref{sec:exp}, one fixed pair $(\alpha_1,\beta_1)=(1,1\times10^{-10})$ works
fairly well. In principle, it is plausible to estimate them from the data simultaneously
with other parameters in order to adaptively accommodate noise features, especially for
large data sets \cite{GelmanCarlinSternRubin:2004}. However, there are no conjugate
priors compatible with the adopted variational framework \cite{Tipping:2005}. Therefore,
one possible way is to maximize the divergence with respect to $(\alpha_1,\beta_1)$. It
is easy to find that in \eqref{eqn:kld}, the only term relates to $\alpha_1$ and
$\beta_1$ is given by
\begin{equation*}
   n\alpha_1\ln \beta_1 - n \ln\Gamma(\alpha_1)+ (\alpha_1-1)\sum_{i=1}^nE_{q^\ast(w_i)}[\ln w_i]-\beta_1\sum_{i=1}^nE_{q^\ast(w_i)}[w_i].
\end{equation*}
Taking derivatives with respect to $\alpha_1$ and $\beta_1$, we arrive at
\begin{equation*}
  \begin{aligned}
    \ln \beta_1 - \psi(\alpha_1) + \frac{1}{n}\sum_{i=1}^nE_{q^\ast(w_i)}[\ln w_i]=0,\\
    \frac{\alpha_1}{\beta_1}-\frac{1}{n}\sum_{i=1}^nE_{q^\ast(w_i)}[w_i]=0,
  \end{aligned}
\end{equation*}
where $\psi(s)=\frac{\partial}{\partial s}\ln\Gamma(s)$ denotes the digamma function. The
solution to the system is not available in closed form, but upon eliminating the variable
$\beta_1$, it can be solved efficiently by the Newton-Raphson method, see e.g.,
\cite[Sect. 3.1]{ChoiWette:1969}.

Finally, we briefly mention the extension to nonlinear problems via recursive
linearization \cite{ChappellGrovesWhitcherWoolrich:2009,JinZou:2010jcp}. The main idea is
to approximate the nonlinear model $\mathbf{K}(\mathbf{u})$ by its first-order Taylor
expansion $\widetilde{\mathbf{K}}(\mathbf{u})$ around the mode $\tilde{\mathbf{u}}$ of an
approximate posterior, i.e.,
\begin{equation*}
  \widetilde{\mathbf{K}}(\mathbf{u})=\mathbf{K}(\tilde{\mathbf{u}})+\mathbf{J}(\mathbf{u}-\tilde{\mathbf{u}}),
\end{equation*}
where $\mathbf{J}=\nabla_\mathbf{u}\mathbf{K(\tilde{\mathbf{u}})}$ is the Jacobian of the
model $\mathbf{K}$ with respect to $\mathbf{u}$. With this linearized model
$\widetilde{\mathbf{K}}(\mathbf{u})$ in place of $\mathbf{K}(\mathbf{u})$, Algorithm
\ref{alg:va} might be employed to deliver an approximation. The mode of the this
newly-derived approximation is then taken for (hopefully) more accurately capturing the
nonlinearity of the genuine model $\mathbf{K}(\mathbf{u})$. This procedure is repeated
until a satisfactory solution is achieved, which gives rise to Algorithm
\ref{alg:nonlin}. In the inner loop, the variational approximation needs not be carried
out very accurately. As to the stopping criterion at Step 7, there are several choices,
e.g., based on the relative change of the inverse solution $\mathbf{u}$.

\begin{algorithm}
\caption{Variational approximation for nonlinear models $\mathbf{K}$}\label{alg:nonlin}
\begin{algorithmic}[1]
 \STATE Given initial guess $q^0(\mathbf{u})$, $q^0(\mathbf{w})$ and $q^0(\lambda)$,
        and set $k=0$.
 \REPEAT
      \STATE Calculate the mode $\tilde{\mathbf{u}}^k$ and the Jacobian
             $\mathbf{J}^k=\nabla_\mathbf{u}\mathbf{K}(\tilde{\mathbf{u}}^k)$.
      \STATE Construct the linearized model, i.e. $\widetilde{\mathbf{K}}(\mathbf{u})
             =\mathbf{K}(\tilde{\mathbf{u}}^k)+\mathbf{J}^k(\mathbf{u}-
             \tilde{\mathbf{u}}^k)$;
      \STATE Find a variational approximation $q^{k+1}(\mathbf{u})q^{k+1}(\mathbf{w})
             q^{k+1}(\lambda)$ using $\widetilde{\mathbf{K}}(\mathbf{u})$ by Algorithm \ref{alg:va};
      \STATE Set $k=k+1$;
 \UNTIL A stopping criterion is satisfied
 \STATE Return $q^k(\mathbf{u})q^k(\mathbf{w})q^k(\lambda)$ as the solution.
\end{algorithmic}
\end{algorithm}

\section{Numerical experiments}\label{sec:exp}

Now we illustrate the proposed method on several benchmark inverse problems in 1d and 2d
heat transfer. These examples are adapted from literature
\cite{SuHewitt:2004,Jin:2008,JinZou:2010jcp,WangZabaras:2004ijhmt,JinZou:2010ima}, and
include both linear and nonlinear models. Throughout, the noisy data $\mathbf{y}$ are
generated as follows
\begin{equation*}
  y_i=\left\{\begin{array}{ll}
     y_i^\dagger, &\mbox{with probability } 1-r\\
     y_i^\dagger+\epsilon\zeta_i, & \mbox{with probability } r
  \end{array}\right.
\end{equation*}
where $\zeta_i$ follow the standard normal distribution, and $(\epsilon,r)$ control the
noise pattern: $r$ is the corruption percentage and $\epsilon=\max_i\{|y_i^\dagger|\}$ is
the corruption magnitude. The matrix $\mathbf{L}$ in the prior $p(\mathbf{u}|\lambda)$ is
taken to be the first-order finite-difference operator, which enforces smoothness on the
sought-for solution $\mathbf{u}$. The parameter pairs $(\alpha_0,\beta_0)$ and
$(\alpha_1,\beta_1)$ are both set to $(1.0,1.0\times10^{-10})$. We have also experimented
with updating $(\alpha_1,\beta_1)$, but it brings little improvement. Hence, we refrain
from presenting the results by adaptively updating $(\alpha_1,\beta_1)$. We shall measure
the accuracy of a solution $\mathbf{u}$ by the relative error
\begin{equation*}
  e=\|\mathbf{u}-\mathbf{u}^\dagger\|_2/\|\mathbf{u}^\dagger\|_2.
\end{equation*}
The algorithm is terminated if the relative change of $\mathbf{u}$ falls below
$tol=1.0\times10^{-5}$. All the computations were performed on a dual core personal
computer with 1.00 GB RAM with MATLAB version 7.0.1 .

\subsection{Cauchy problem}
\begin{figure}
  \begin{tabular}{cc}
     \includegraphics[width=7cm]{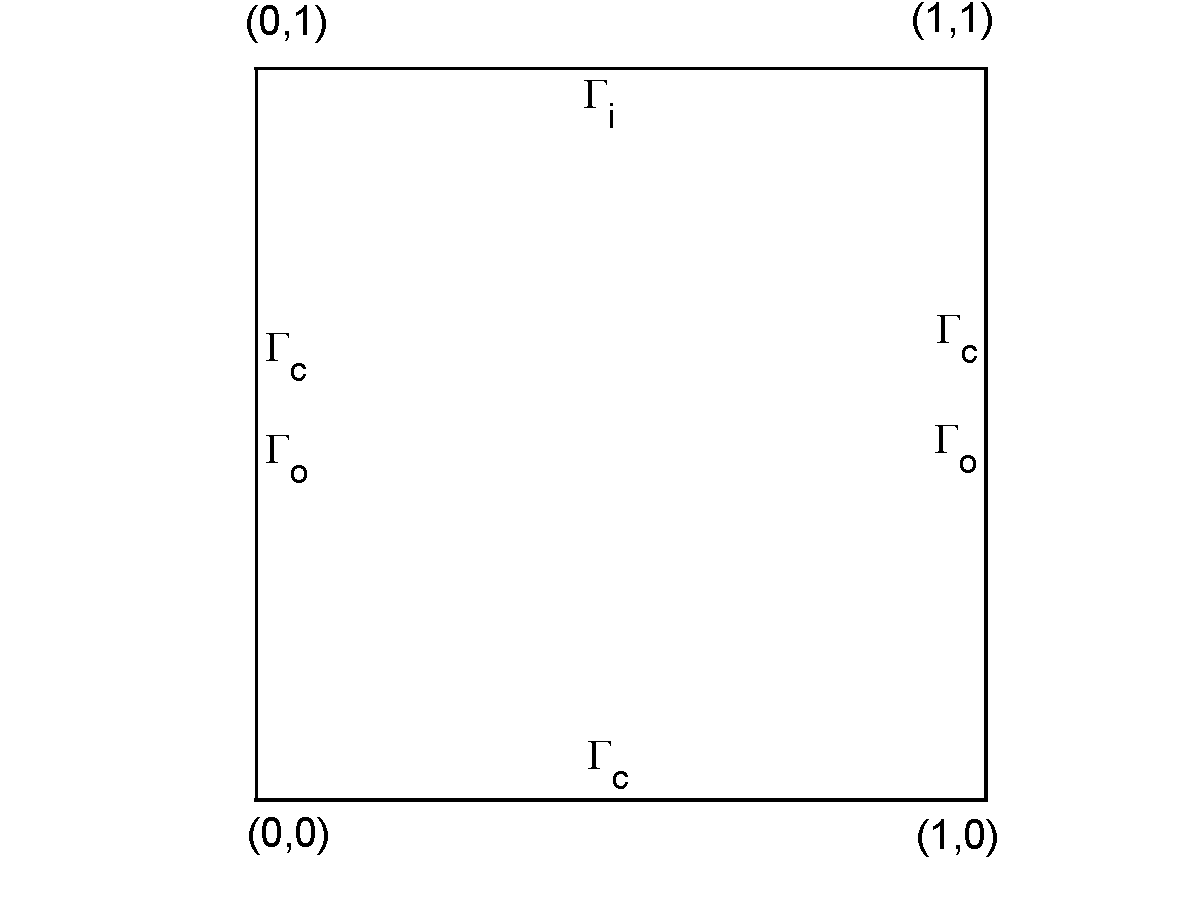} & \includegraphics[width=7cm]{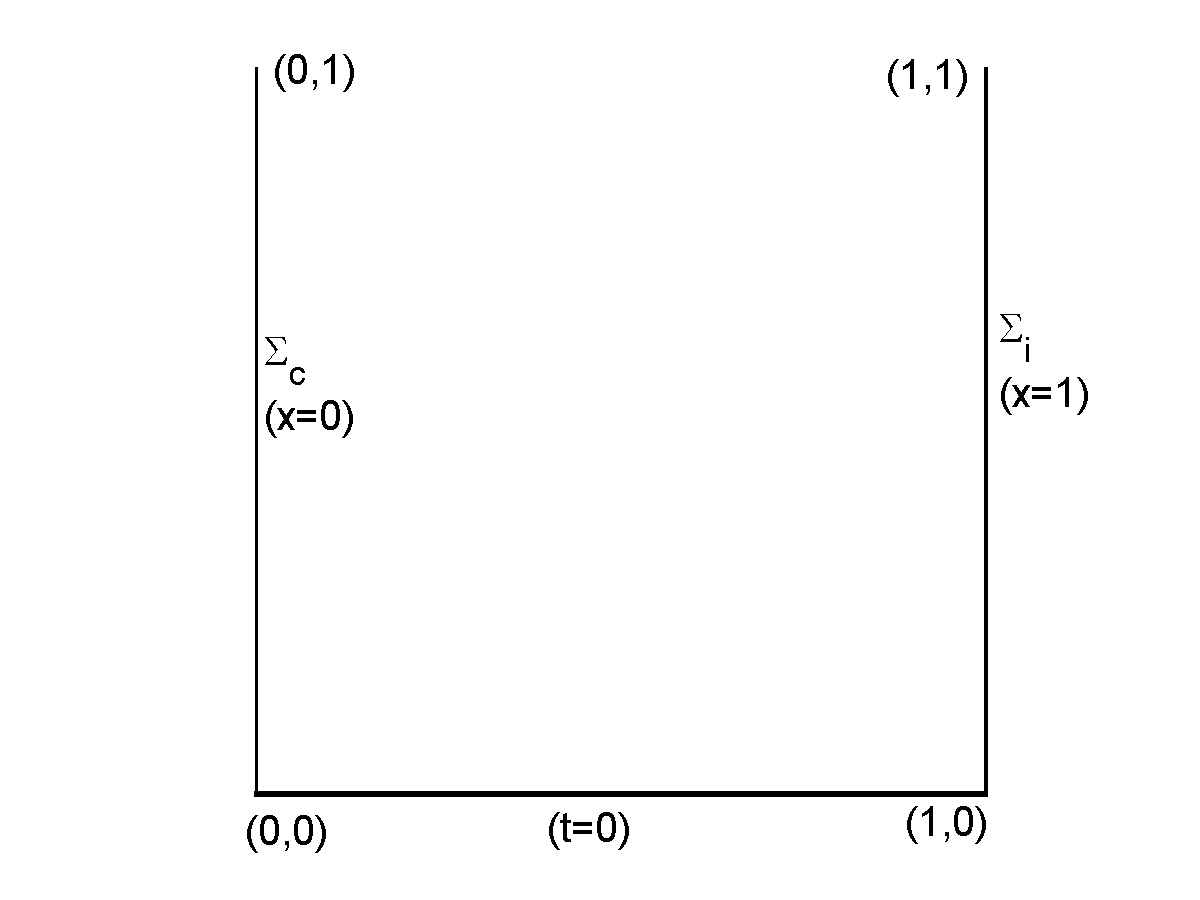}\\
       (a) & (b)
  \end{tabular}
  \caption{The domain and boundaries for Example 1 (a) and Example 3 (b).}\label{fig:domain}
\end{figure}

This example is taken from \cite[Sect. 5.2.1]{JinZou:2010jcp}. Here we consider the
Cauchy problem for steady state heat conduction. Let $\Omega$ be the unit square
$(0,1)^2$ with its boundary $\Gamma$ divided into two disjoint parts, i.e.,
$\Gamma_i=[0,1]\times\{1\}$ and $\Gamma_c = \Gamma\backslash\Gamma_i$. The steady-state
heat conduction is described by
\begin{equation*}
  -\Delta y = 0\quad \mbox{in } \Omega.
\end{equation*}
It is subjected to the boundary conditions
\begin{equation*}
\frac{\partial y}{\partial n}=g \quad\mbox{on}\quad \Gamma_c \quad \quad \mbox{and}\quad \quad
y=u\quad \mbox{on}\quad \Gamma_i,
\end{equation*}
where $n$ denotes the unit outward normal. The linear operator $K$ maps $u$ (with $g=0$)
to $y$ restricted to the segments $\Gamma_o=\{0,1\}\times(0, 1)\subset\Gamma_c$. We refer
to Fig. \ref{fig:domain}(a) for a schematic plot of the domain and its boundaries. The
inverse problem seeks the unknown $u$ from noisy data $y$. It arises, e.g., in the study
of re-entrant space shuttles \cite{Beck:1985} and electro-cardiography
\cite{ColliGuerri:1985}. For the inversion, the solution $y$ is given by $\sin\pi
x_1e^{\pi x_2}+x_1+x_2$, from which both $g$ and $u$ can be evaluated directly. The
operator $K$ is discretized using piecewise linear finite element with $3200$ triangular
elements, see \cite{JinZou:2010jcp,JinZou:2010ima} for details. The number of
measurements $\mathbf{y}$ is $80$, and the unknown $\mathbf{u}$ is of dimension $41$.

The numerical results for the example with various levels of noises are shown in Table
\ref{tab:exam1}, where $e$ refers to the relative error. A first observation is that the
corruption percentage $r$ plays an important role in the error $e$, and there is a sudden
loss of the accuracy $e$ as $r$ increases from $0.4$ to $0.5$. Nonetheless, the
reconstruction remains very accurate for $r$ up to $0.9$.

\begin{table}
\centering
  \caption{Numerical results for Example 1 with various noise levels.}\label{tab:exam1}
  \begin{tabular}{c|ccccccccc}
   \toprule
        $r$ & $0.1$   & $0.2$   & $0.3$   & $0.4$   & $0.5$   & $0.6$   & $0.7$   & $0.8$   & $0.9$  \\
   \midrule
   $\lambda$& 8.56e-1 & 8.51e-1 & 8.49e-1 & 8.49e-1 & 8.36e-1 & 8.28e-1 & 8.28e-1 & 8.28e-1 & 8.27e-1\\
   $e$      & 2.33e-4 & 3.67e-4 & 3.65e-4 & 3.67e-4 & 1.59e-3 & 2.49e-3 & 2.49e-3 & 2.49e-3 & 2.53e-3\\
  \bottomrule
  \end{tabular}
\end{table}

\begin{figure}
  \centering
  \begin{tabular}{cc}
  \includegraphics[width=7cm]{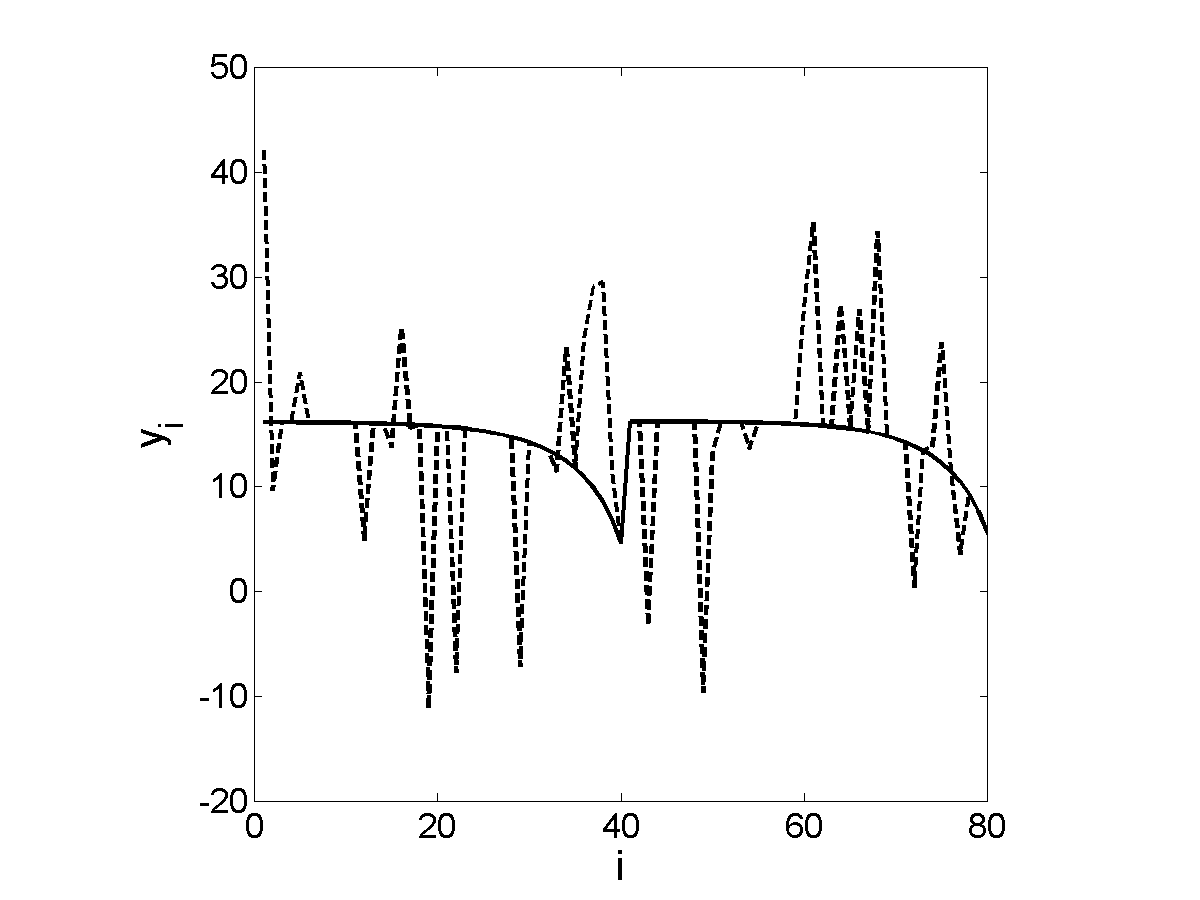} & \includegraphics[width=7cm]{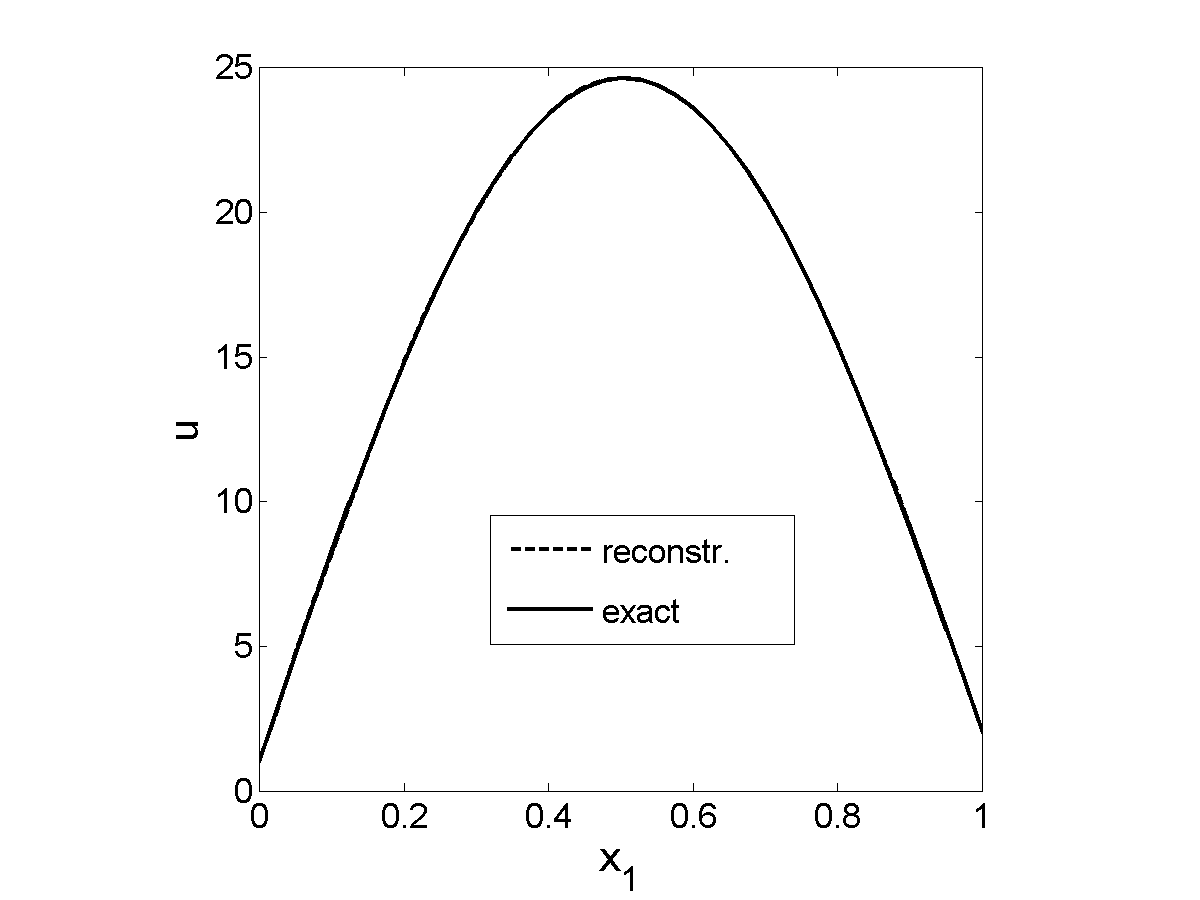}\\
  (a) noisy data of level $r=0.5$ & (b) mean by $t$ model \\
  \includegraphics[width=7cm]{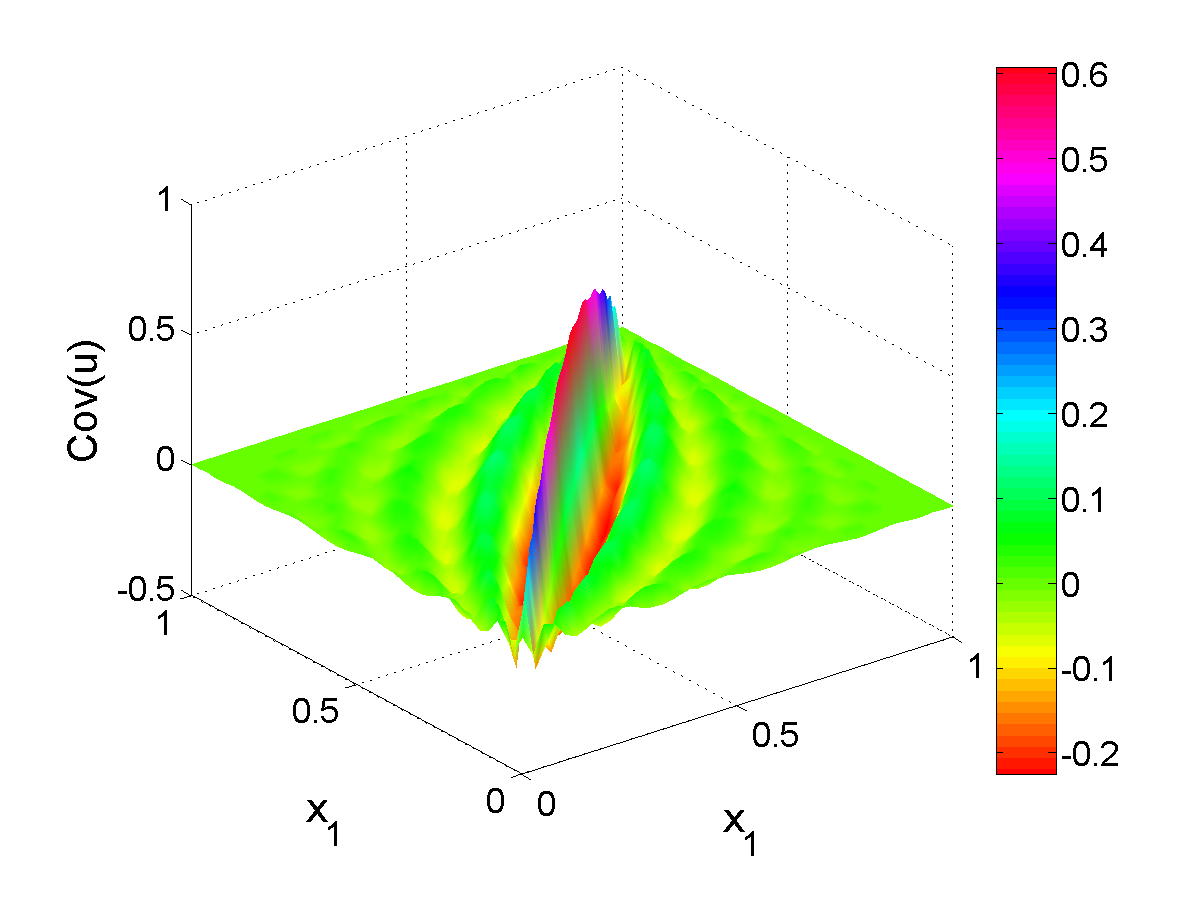} &  \includegraphics[width=7cm]{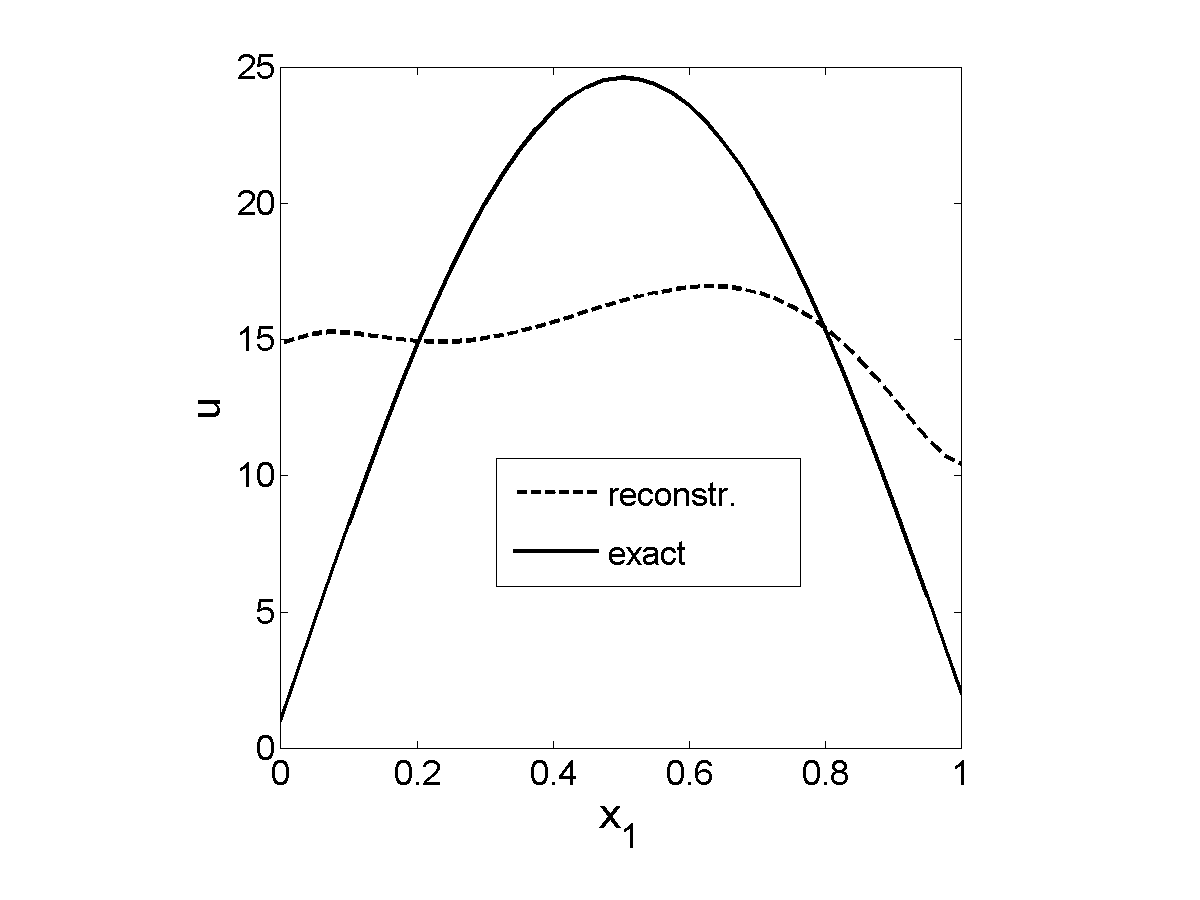}\\
  (c) covariance by $t$ model & (d) mean by Gaussian model\\
  \includegraphics[width=7cm]{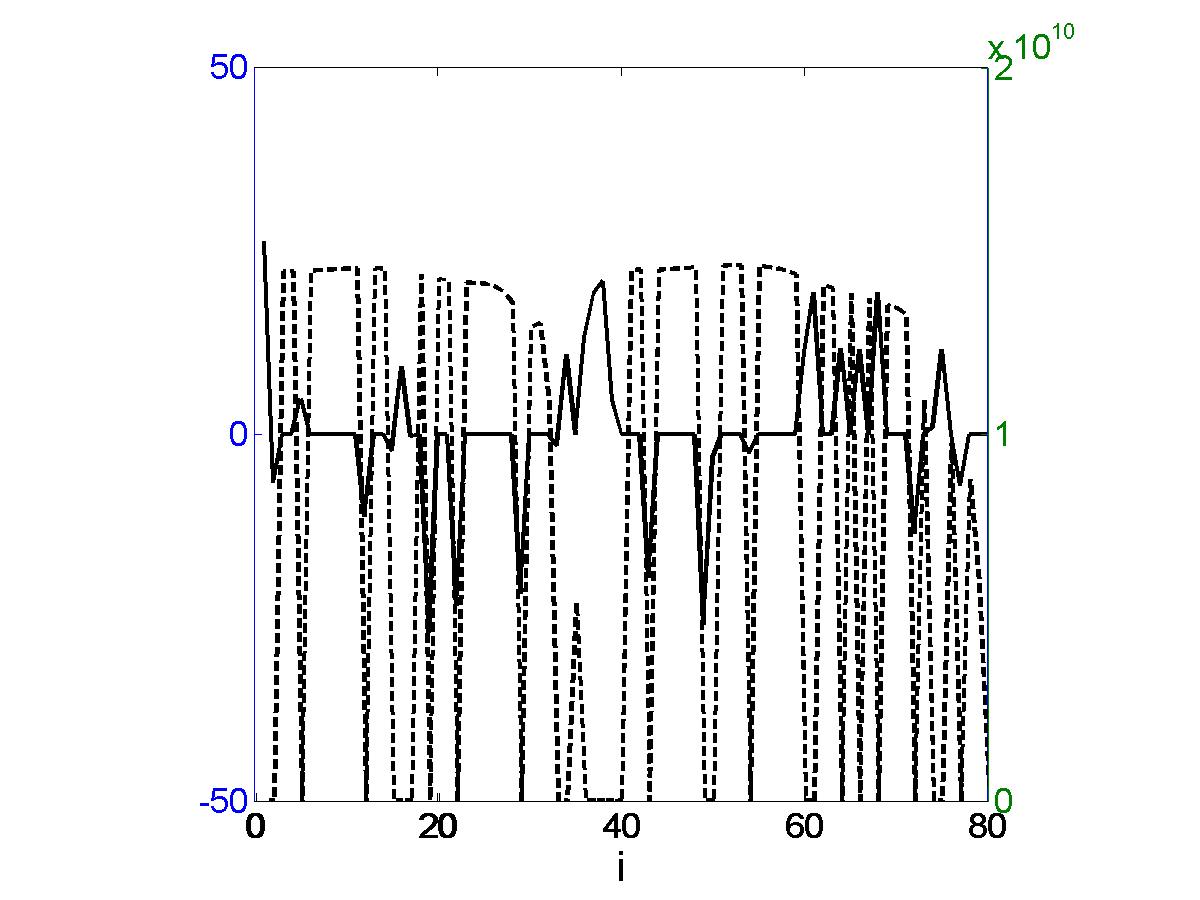} &   \includegraphics[width=7cm]{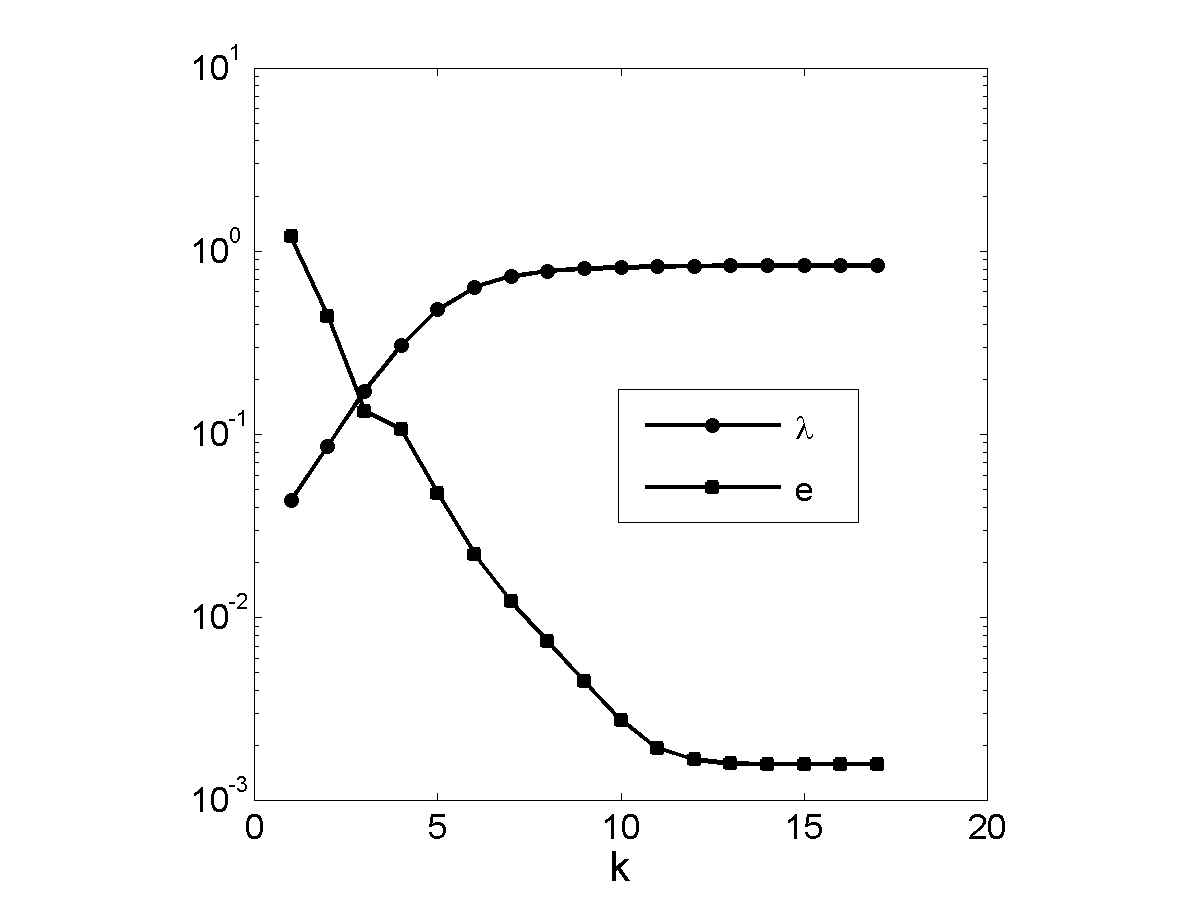}\\
  (e) weight $\mathbf{w}$ v.s. noise $\boldsymbol\zeta$ & (f) Convergence of Alg. \ref{alg:va}
  \end{tabular}
  \caption{Numerical results for Example 1 with $r=0.5$ noise. In (e), the solid and dashed lines
  refer to the noise $\boldsymbol\zeta$ and the weight $\mathbf{w}$, respectively.}\label{fig:exam1sol}
\end{figure}

A typical realization of the noisy data of level $r=0.5$ is shown in Fig.
\ref{fig:exam1sol}(a). We observe that some data points deviate significantly from the
exact ones, and are completely erroneous. The solutions (mean of the Bayesian solution)
by the proposed approach ($t$ model) and the Gaussian model is shown in Figs.
\ref{fig:exam1sol}(b) and \ref{fig:exam1sol}(d), respectively, where $x_1$ is the first
spatial coordinate. Here, for the Gaussian model, the regularization parameter $\eta$ is
manually tuned, i.e. $\eta=4.64$, so as to yield a solution with the smallest possible
error. The solution by the $t$ model is in excellent agreement with the exact one, while
that by the standard approach is completely off the track. This shows clearly the
robustness of the $t$ model. The covariance, see Fig. \ref{fig:exam1sol}(c), can be used
for quantifying the uncertainty of a specific solution. The covariance is relatively
smooth, and decays quickly away from neighboring sites. The weight $\mathbf{w}$ admits
nice interpretations. We plot in Fig. \ref{fig:exam1sol}(e) the noise (solid line, value
in blue) and the weight (dashed line, value in green). There is a one-to-one
correspondence of the noise sites and the sites of the weight with small values. Thus the
weight $\mathbf{w}$ can accurately detect the locations of the noises and effectively
prunes out the noises from the inversion procedure simultaneously. The convergence of
Algorithm \ref{alg:va} is very steady and fast, see Fig. \ref{fig:exam1sol}(f), and it is
reached within about ten iterations. The convergence of the scalar $\lambda$ seems
monotonic.

\subsection{Flux reconstruction}
\begin{figure}
  \centering
  \begin{tabular}{cc}
    \includegraphics[width=7cm]{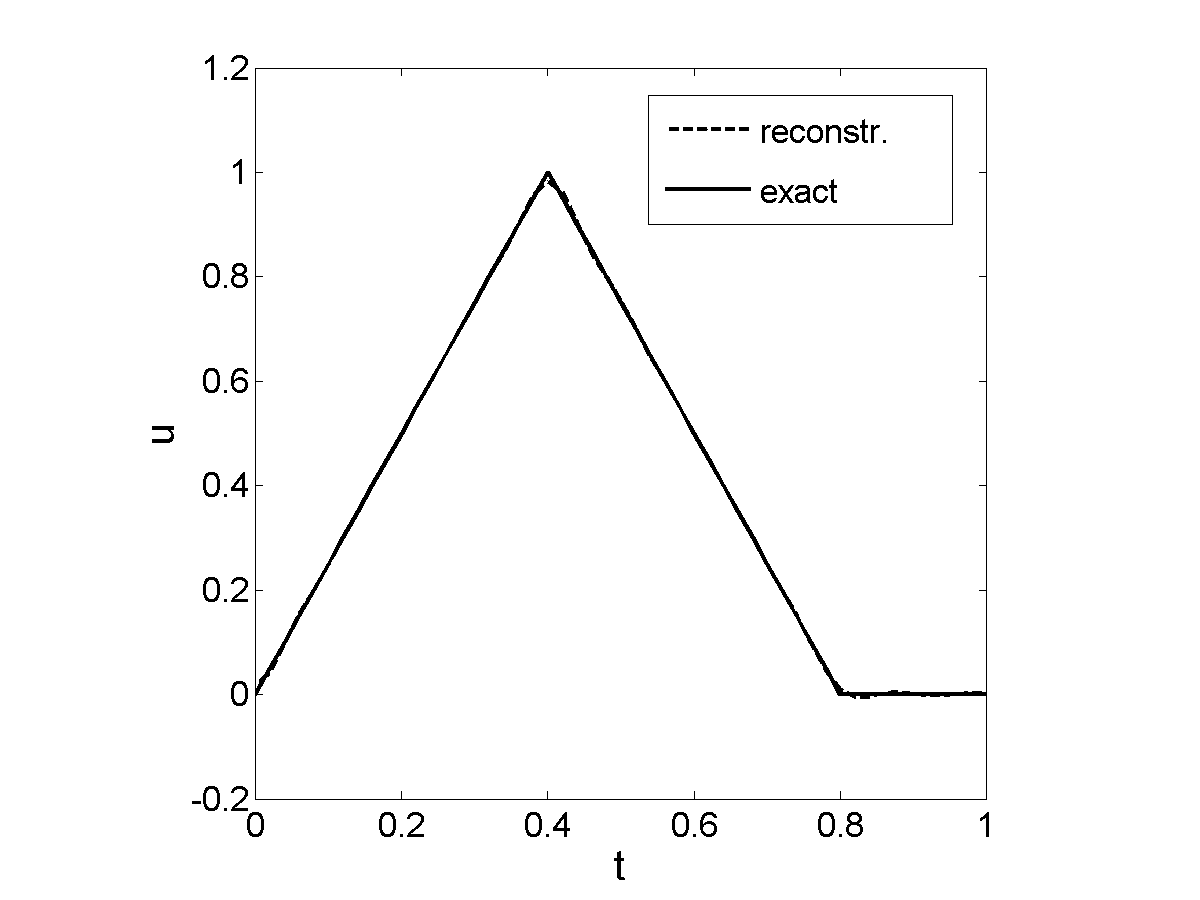} &   \includegraphics[width=7cm]{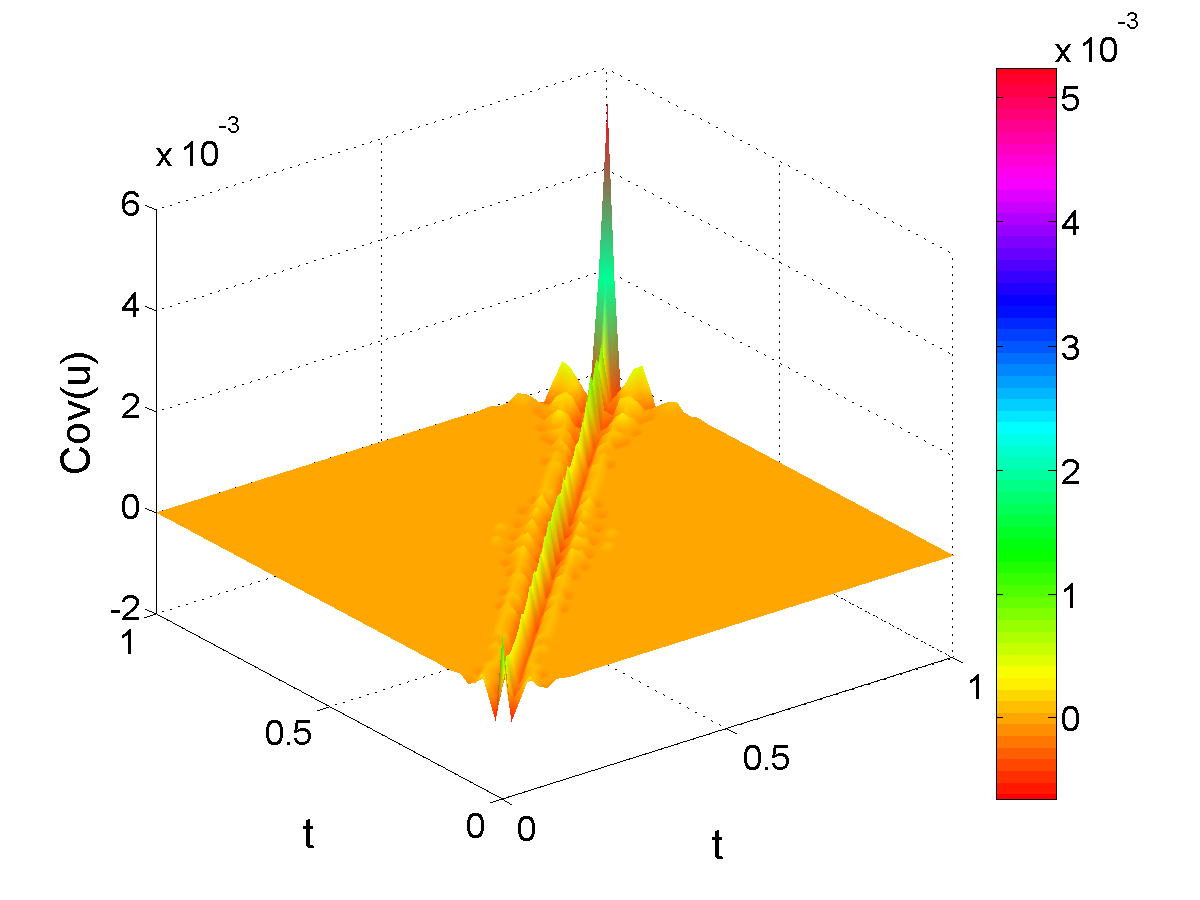}\\
    (a) mean by $t$ model & (b) covariance by $t$ model\\
    \includegraphics[width=7cm]{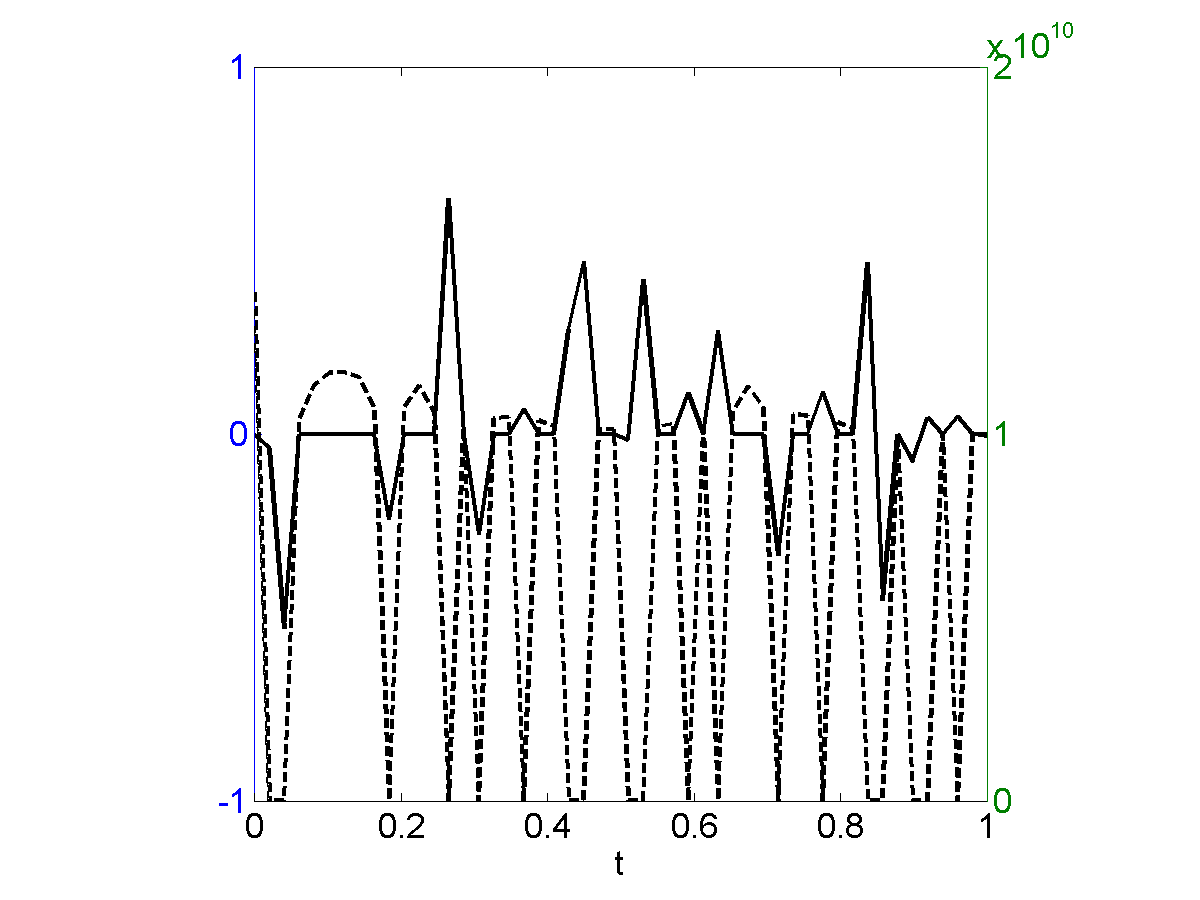} &  \includegraphics[width=7cm]{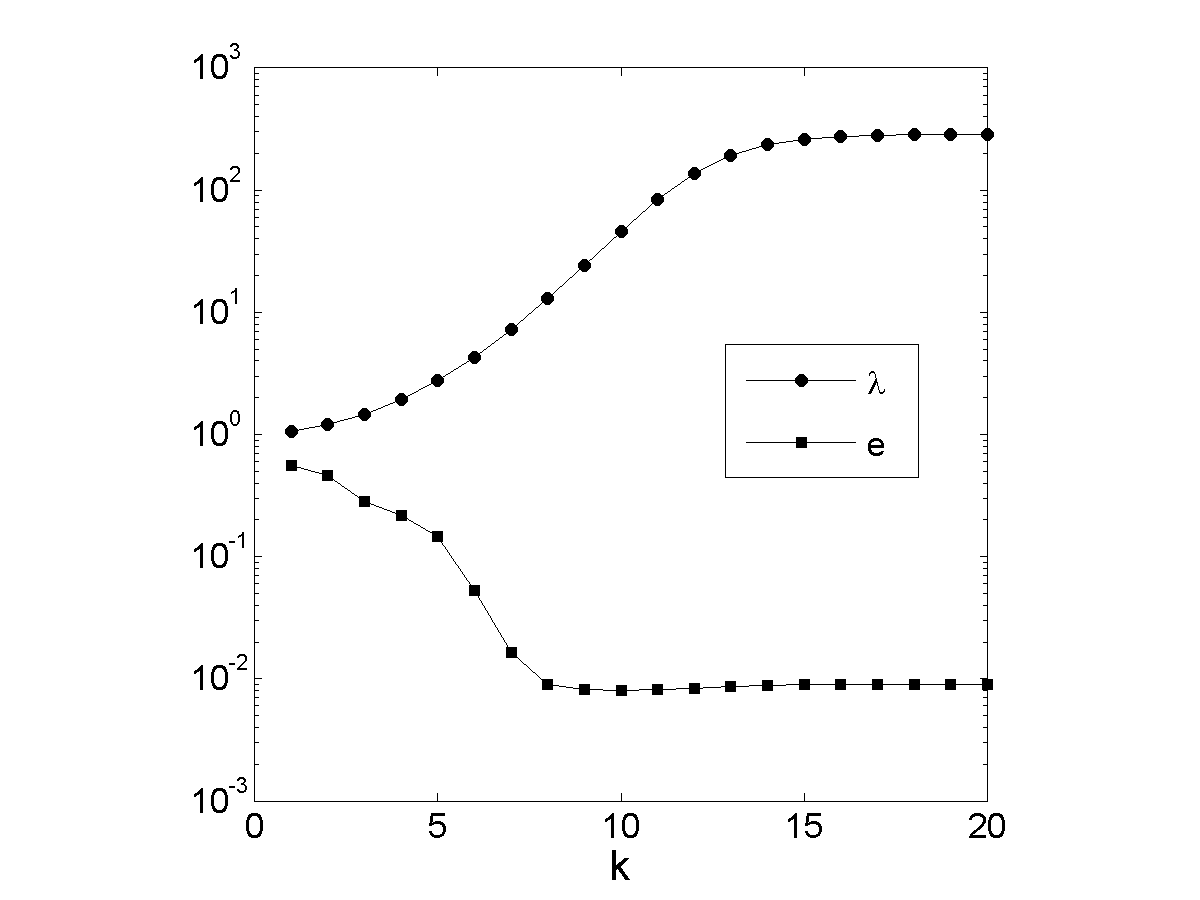}\\
    (c) weight $\mathbf{w}$ v.s. noise $\boldsymbol\zeta$ & (d) convergence of Alg. \ref{alg:va}
  \end{tabular}
  \caption{Numerical results for Example 2 with $r=0.5$ noise. In (c), the solid and dashed lines
  refer to the noise $\boldsymbol\zeta$ and the weight $\mathbf{w}$, respectively.}\label{fig:exam2sol}
\end{figure}

This example is taken from \cite[Sect. 5.1]{WangZabaras:2004ijhmt}. Here, we consider 1d
transient heat transfer. Let $\Omega$ be the interval $(0,1)$, and the time interval be
$[0,1]$. The 1d transient heat conduction is described by
\begin{equation*}
  \frac{\partial y}{\partial t}=\Delta y,
\end{equation*}
with a zero initial condition and the following boundary conditions
\begin{equation*}
 \frac{\partial y}{\partial n}=g(t) \quad \mbox{on} \quad \Sigma_c \quad\quad\mbox{and}\quad\quad
 \frac{\partial y}{\partial n}=u(t) \quad\mbox{on}\quad\Sigma_i,
\end{equation*}
where the boundaries $\Sigma_c=\{x=0\}\times[0,T]$ and $\Sigma_i=\{x=1\}\times[0,T]$. The
operator $K$ maps the flux $u$ (with $g=0$) to $y$ restricted to $\Sigma_c$. The inverse
problem is to recover the flux $u$ from noisy data $y$. For the inversion, we take
$g(t)=0$ and a hat shaped flux $u$, see Fig. \ref{fig:exam2sol}(b) for its profile. The
spatial and temporal intervals are discretized into $101$ and $201$ uniform grids,
respectively. The operator $K$ is discretized with piecewise linear finite elements in
space and backward finite-difference in time. The number of measurements $\mathbf{y}$ is
$50$, and the unknown $\mathbf{u}$ is on a coarse mesh and of size $51$.

\begin{table}
  \centering
  \caption{Numerical results for Example 2 with various noise levels.}\label{tab:exam2}
  \begin{tabular}{c|cccccccc}
  \toprule
        $r$& $0.1$   & $0.2$   & $0.3$   & $0.4$   & $0.5$   & $0.6$   & $0.7$   & $0.8$  \\
  \midrule
  $\lambda$& 3.29e2  & 3.20e2  & 3.04e2  & 2.97e2  & 2.86e2  & 2.71e2  & 2.55e2  & 2.44e2\\
        $e$& 5.51e-3 & 6.73e-3 & 8.17e-3 & 8.28e-3 & 9.06e-3 & 9.41e-3 & 1.79e-2 & 1.79e-2\\
  \bottomrule
  \end{tabular}
\end{table}

The numerical results for Example 2 are shown in Table \ref{tab:exam2}. The $\lambda$
value is independent of the corruption percentage $r$, and the accuracy $e$ only
deteriorates very mildly with the increase of the corruption percentage $r$ from $0.1$ to
$0.8$. The solution for a typical realization of noisy data of level $r=0.5$ is shown in
Fig. \ref{fig:exam2sol}(a), where $t$ is the temporal coordinate. It agrees excellently
with the true solution, except small errors around the corner. The variance at the end
points, especially around $t=1$, is much more pronounced than that in the interior, see
Fig. \ref{fig:exam2sol}(b). This might be related to the causality nature of heat
problems. The weight $\mathbf{w}$ detects noise sites accurately and meanwhile eliminates
them from the inversion by putting very small weight, and the convergence of the
algorithm is steady and fast, c.f. Figs. \ref{fig:exam2sol}(c) and (d), respectively.

\subsection{Stationary Robin inverse problem}

This example is adapted from \cite[Sect. 5.2.4]{JinZou:2010jcp} \cite{Jin:2008}, to
illustrate the approach for nonlinear problems. Let $\Omega$ be the unit square $(0,1)^2$
with its boundary $\Gamma$ divided into two disjoint parts, i.e.,
$\Gamma_i=[0,1]\times\{1\}$ and $\Gamma_c = \Gamma\backslash\Gamma_i$. The steady-state
heat conduction is described by
\begin{equation*}
  -\Delta y = 0\quad \mbox{in } \Omega.
\end{equation*}
It is equipped with the following boundary conditions
\begin{equation*}
  \frac{\partial y}{\partial n} = g(x)\quad \mbox{on}\quad \Gamma_c \quad\quad\mbox{and}\quad\quad
  \frac{\partial y}{\partial n} + u y =0\quad \mbox{on}\quad \Gamma_i,
\end{equation*}
where $u$ is the heat transfer coefficient. The operator $K$ maps the coefficient $u$ to
$y$ restricted to $\Gamma_c$. The inverse problem is to reconstruct the unknown $u$ from
noisy data $y$. It arises in corrosion detection \cite{Inglese:1997,JinZou:2010ima} and
analysis of quenching process \cite{OsmanBeck:1988}. For the inversion, the flux $g$ is
set to $1$, and the true coefficient $u$ is given by $1+\sin(\pi x_1)$. The operator $K$
is discretized using piecewise linear finite element with $3200$ triangular elements. The
number of measurements $\mathbf{y}$ is $120$, and the unknown $\mathbf{u}$ is of
dimension $41$.

\begin{figure}
  \centering
  \begin{tabular}{cc}
    \includegraphics[width=7cm]{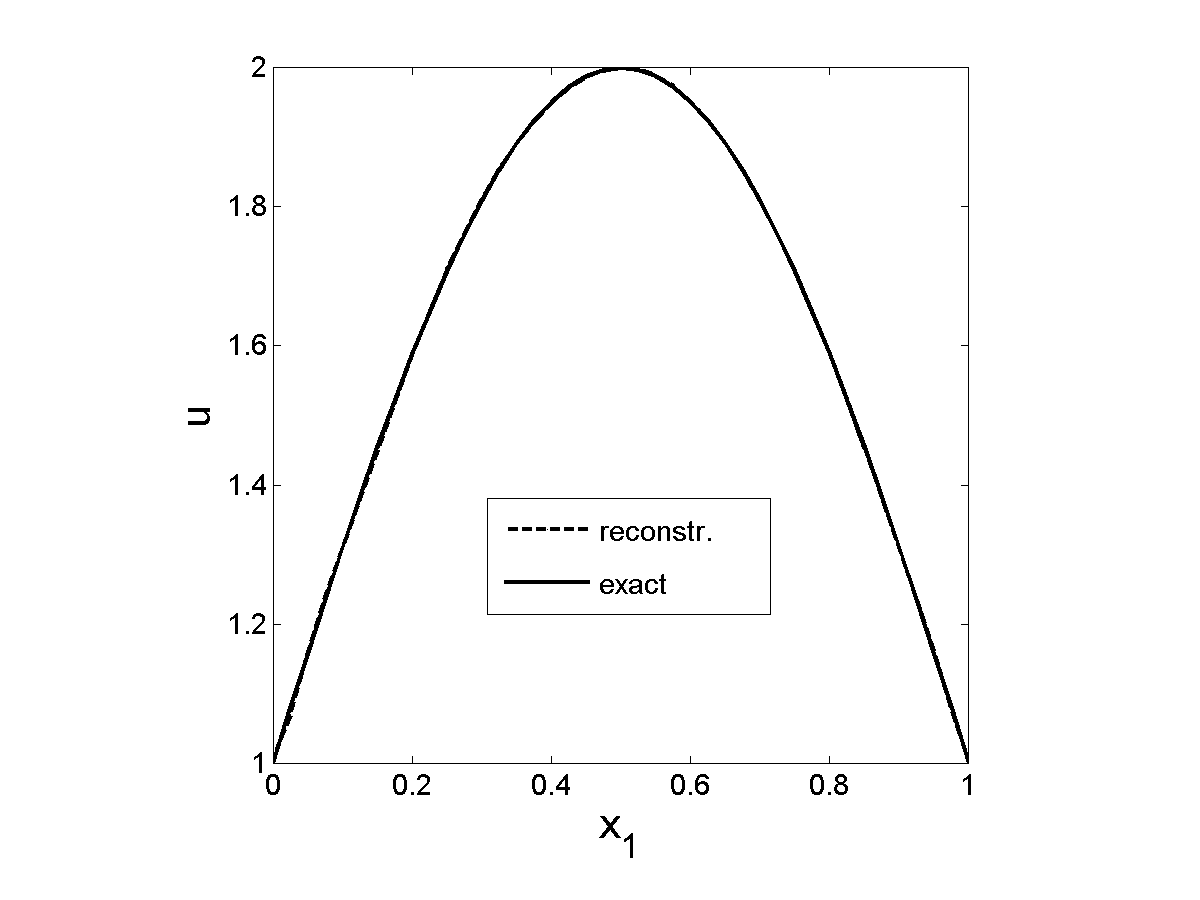} & \includegraphics[width=7cm]{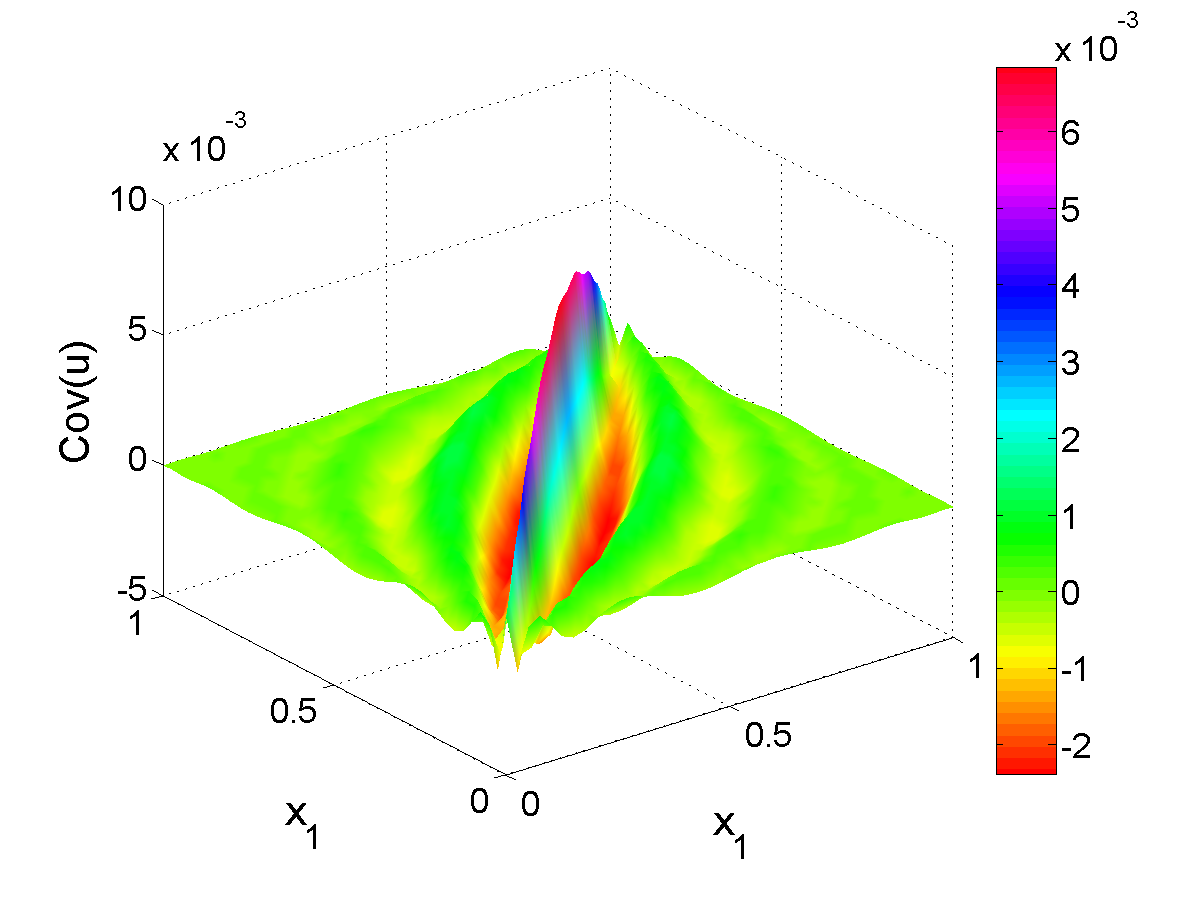}\\
    (a) mean by $t$ model & (b) covariance by $t$ model \\
    \includegraphics[width=7cm]{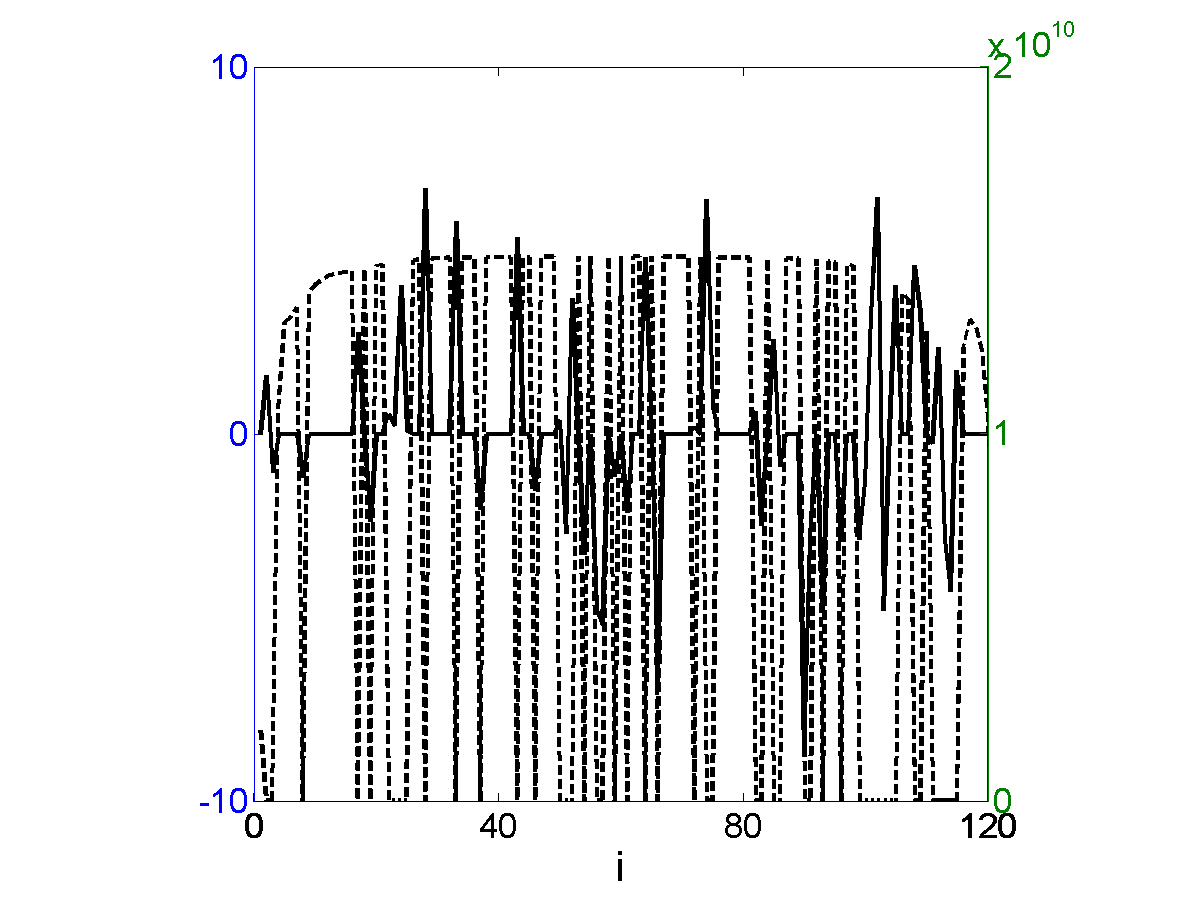} & \includegraphics[width=7cm]{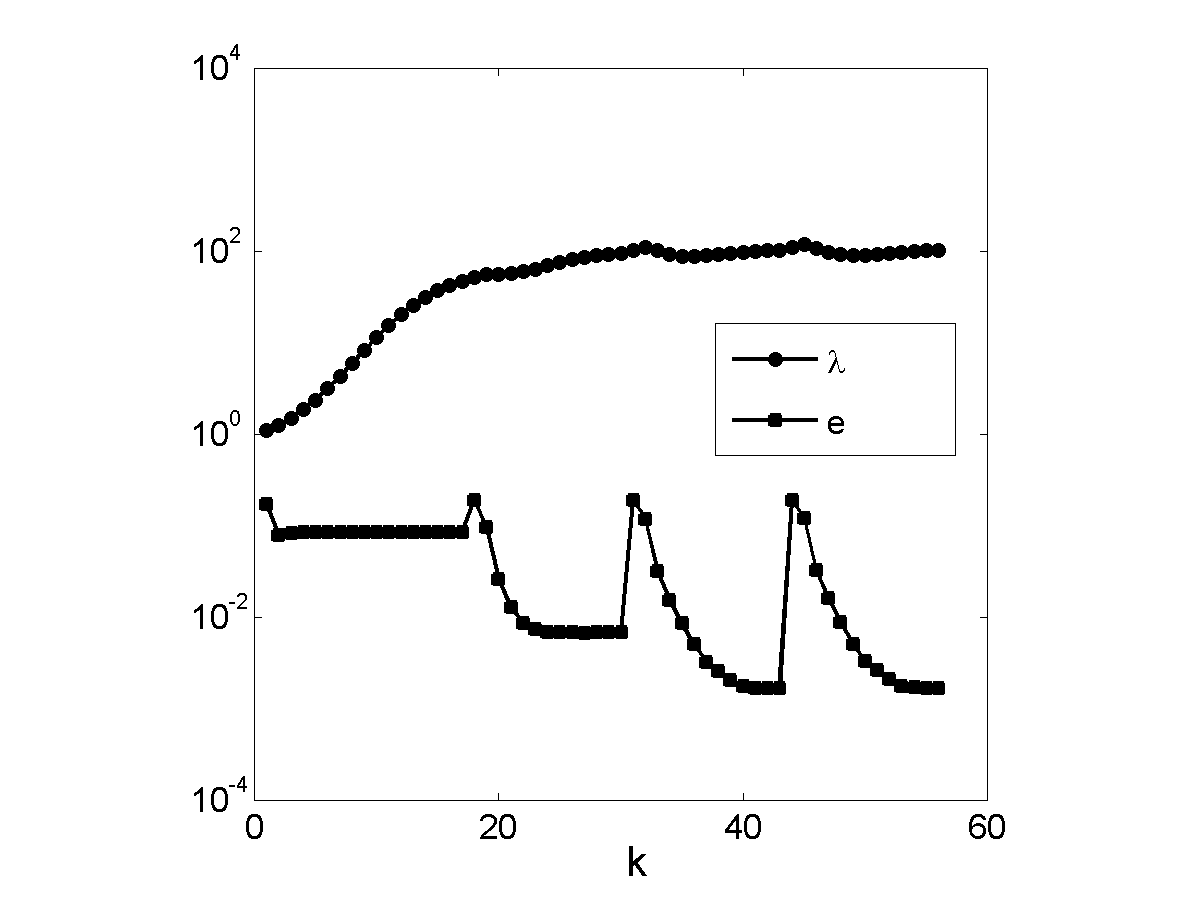}\\
    (c) weight $\mathbf{w}$ v.s. noise $\boldsymbol\zeta$ & (d) convergence of Alg. \ref{alg:nonlin}
  \end{tabular}
  \caption{Numerical results for Example 3 with $r=0.5$ noise.  In (c), the solid and dashed lines
  refer to the noise $\boldsymbol\zeta$ and the weight $\mathbf{w}$, respectively.}\label{fig:exam3sol}
\end{figure}

\begin{table}
  \centering
  \caption{Numerical results for Example 3 with various noise levels.}\label{tab:exam3}
  \begin{tabular}{c|ccccccccc}
  \toprule
        $r$  & $0.1$   & $0.2$   & $0.3$   & $0.4$   & $0.5$   & $0.6$   & $0.7$   & $0.8$   & $0.9$   \\
  \midrule
    $\lambda$& 1.07e2  & 1.04e2  & 1.03e2  & 1.03e2  & 1.03e2  & 1.03e2  & 1.02e2  & 1.01e2  & 9.78e1  \\
         $e$ & 1.30e-3 & 1.72e-3 & 1.71e-3 & 1.71e-3 & 1.70e-3 & 1.70e-3 & 1.69e-3 & 1.76e-3 & 2.27e-3 \\
  \bottomrule
  \end{tabular}
\end{table}

The numerical results for Example 3 are shown in Table \ref{tab:exam3}. The observations
for the linear models remain valid. The solution for an exemplary noise realization of
level $r=0.5$ is shown in Fig. \ref{fig:exam3sol}(a), which agrees well with the exact
one. The convergence of Algorithm \ref{alg:nonlin} is achieved within four (outer)
iterations, see Fig. \ref{fig:exam3sol}(d). The convergence behavior of the algorithm in
the inner loop is similar to that of linear cases. The plateaus indicate that the
tolerance $tol=1.0\times10^{-5}$ is a bit too conservative, and the first two iterations
may be solved less accurately without sacrificing the accuracy of the final solution.

\subsection{Transient Robin inverse problem}
This last example is adapted from \cite{SuHewitt:2004}. Here we consider again 1d
transient heat transfer. Let $\Omega$ be the spatial interval $(0,1)$, and the time
interval be $[0,1]$. The 1d transient heat conduction is described by
\begin{equation*}
  \frac{\partial y}{\partial t}=\Delta y,
\end{equation*}
with a zero initial condition and the following boundary conditions
\begin{equation*}
   \frac{\partial y}{\partial n}=g(t) \quad \mbox{on}\quad\Sigma_c\quad\quad\mbox{and}\quad\quad
   \frac{\partial y}{\partial n} + uy = 0\quad \mbox{ on }\Sigma_i,
\end{equation*}
where $u(t)$ is a time-dependent heat transfer coefficient, and the boundaries
$\Sigma_c=\{x=0\}\times[0,T]$ and $\Sigma_i=\{x=1\}\times[0,T]$. The operator $K$ maps
the coefficient $u$ to $y$ restricted to $\Sigma_c$. The inverse problem is to estimate
the coefficient $u$ from noisy data $y$ \cite{SuHewitt:2004,JinLu:2011}. For the
inversion, the flux $g$ is set to $1$, and the true coefficient $u=1+\frac{1}{2}
\chi_{[\frac{3}{10},\frac{7}{10}]}$ is discontinuous, where $\chi$ denotes the
characteristic function. The spatial and temporal intervals are both discretized into
$100$ uniform intervals. The operator $K$ is discretized with piecewise linear finite
elements in space and backward finite-difference in time. The number of measurements
$\mathbf{y}$ is $101$, and the the unknown $\mathbf{u}$ is of dimension $101$.

\begin{table}
  \centering
  \caption{Numerical results for Example 4 with various levels.}\label{tab:exam4}
  \begin{tabular}{c|ccccccccc}
  \toprule
    $r$    & $0.1$   & $0.2$   & $0.3$   & $0.4$   & $0.5$   & $0.6$   & $0.7$   & $0.8$   & $0.9$  \\
  \midrule
  $\lambda$& 2.25e2  & 2.23e2  & 2.23e2  & 2.26e2  & 2.27e2  & 2.27e2  & 2.32e2  & 2.34e2  & 2.29e2 \\
    $e$    & 3.94e-2 & 3.97e-2 & 3.98e-2 & 4.05e-2 & 4.11e-2 & 4.14e-2 & 4.25e-2 & 4.31e-2 & 4.32e-2\\
  \bottomrule
  \end{tabular}
\end{table}

\begin{figure}
  \centering
  \begin{tabular}{cc}
  \includegraphics[width=7cm]{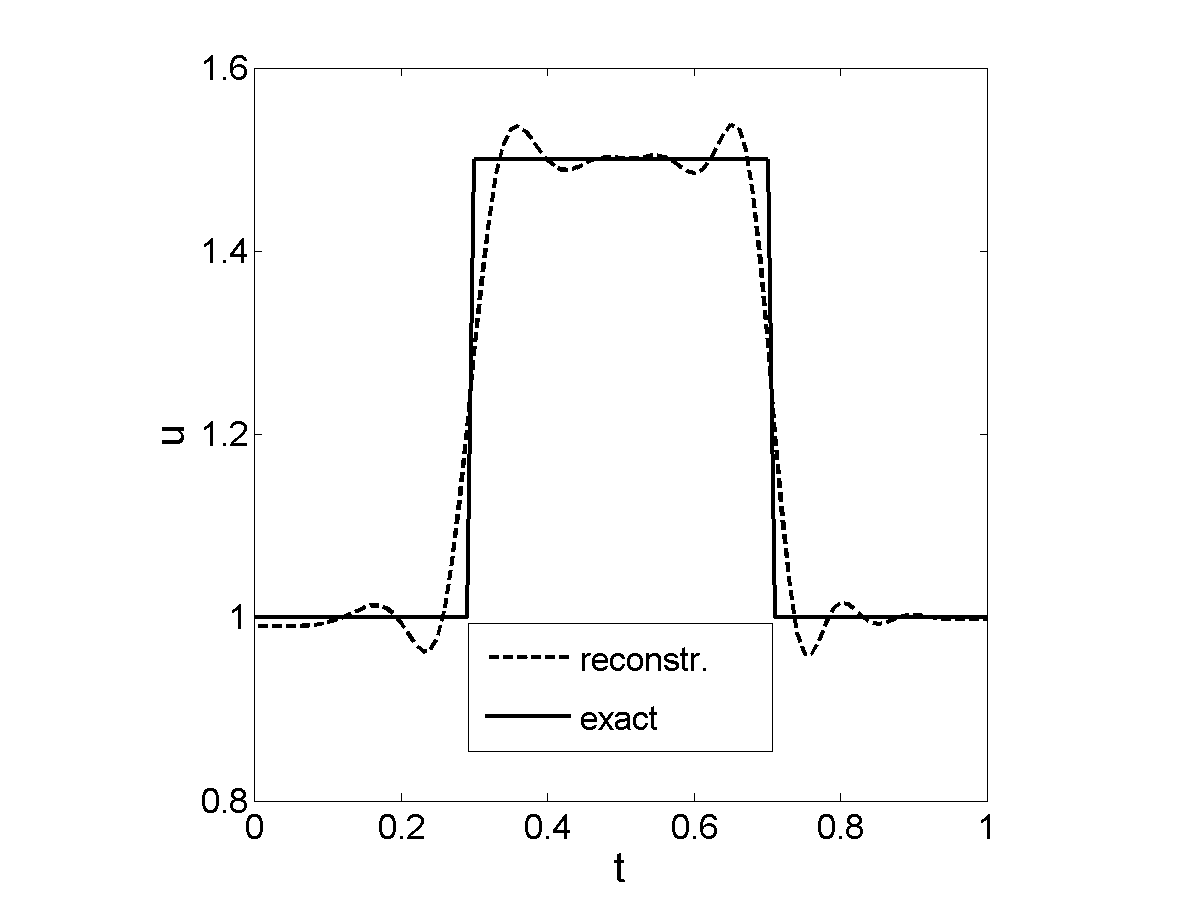} & \includegraphics[width=7cm]{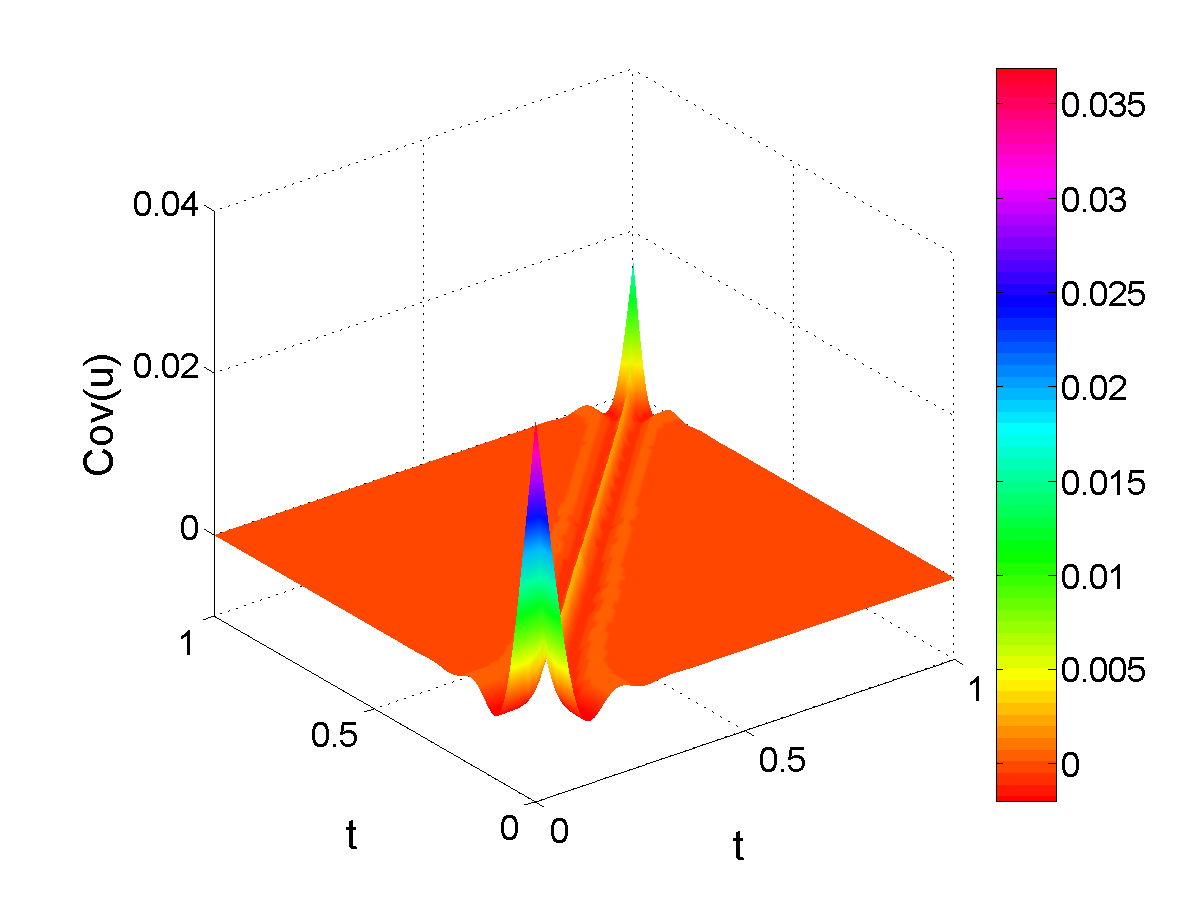}\\
  (a) mean by $t$ model & (b) covariance by $t$ model\\
  \includegraphics[width=7cm]{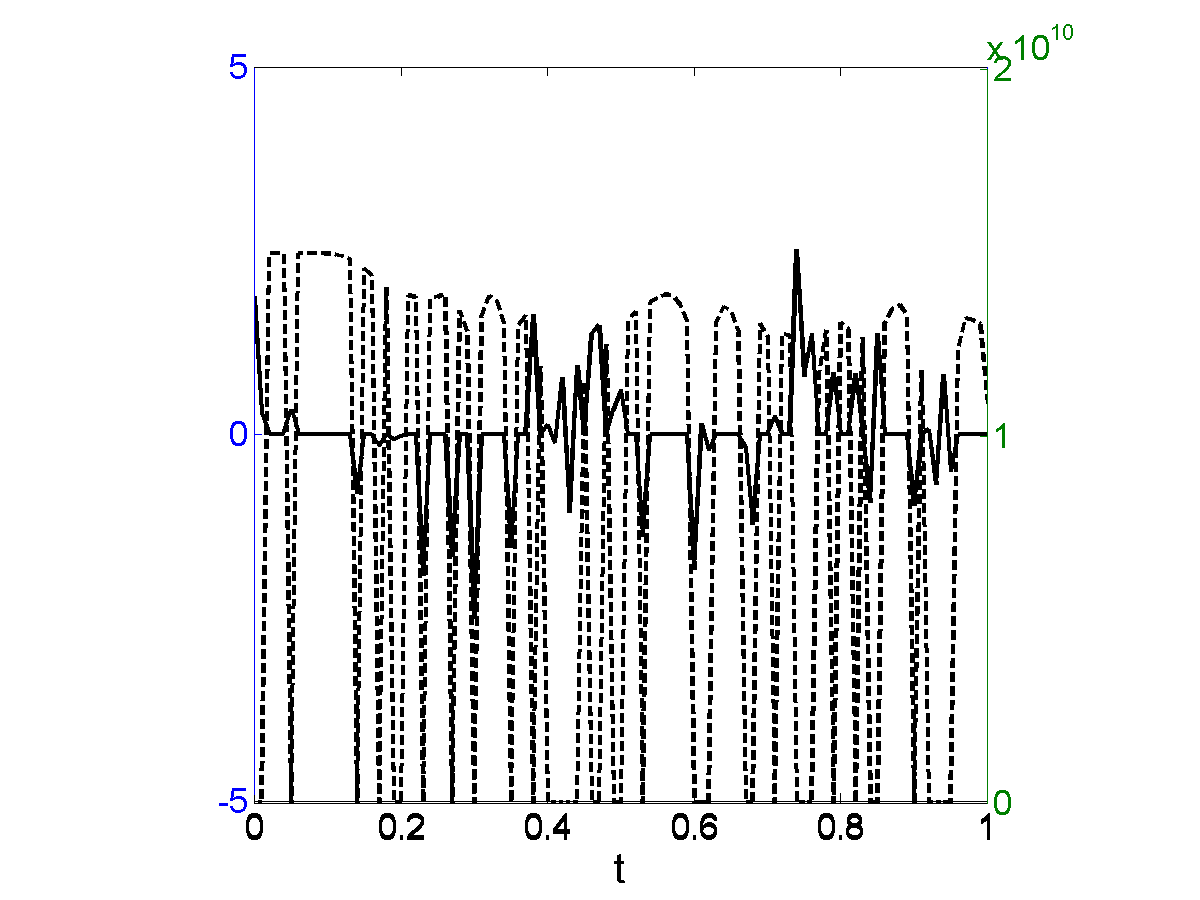} & \includegraphics[width=7cm]{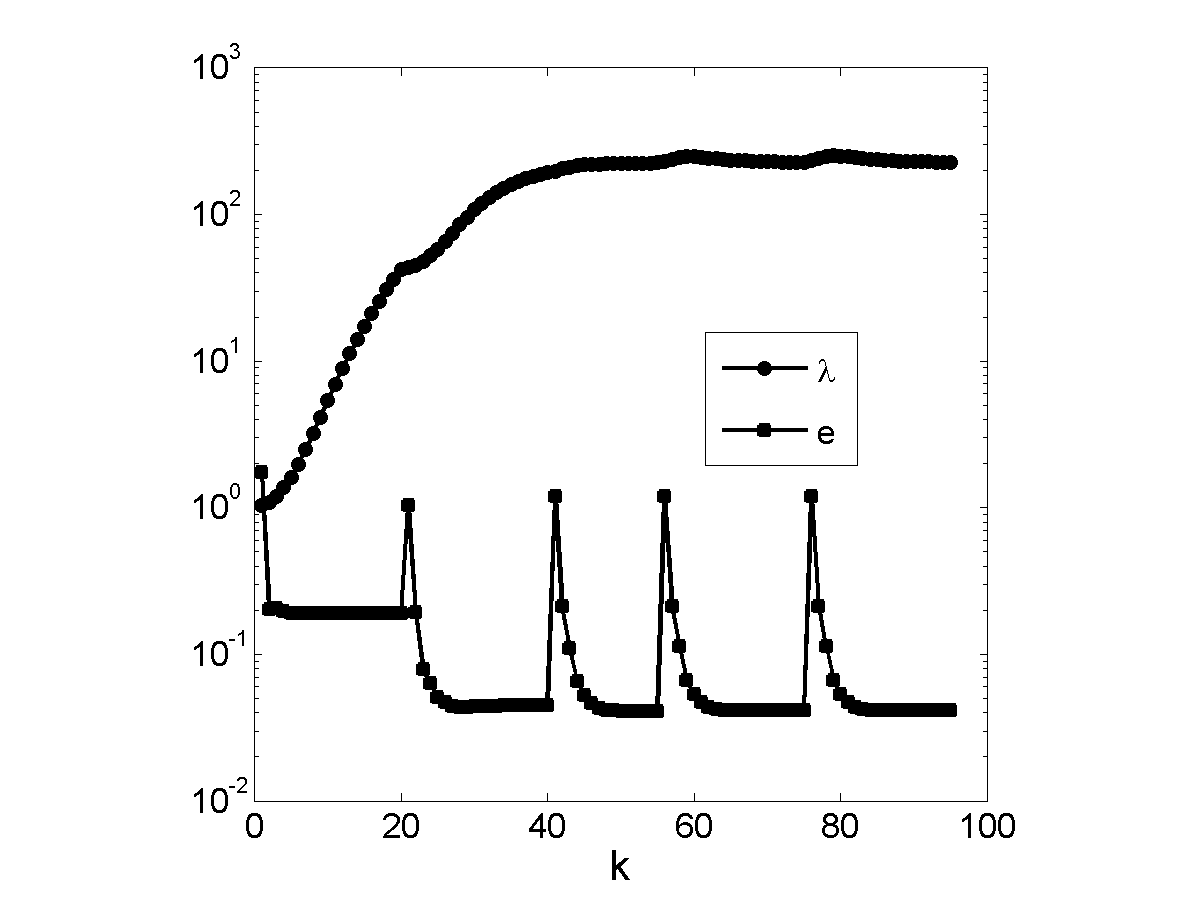}\\
  (c) weight $\mathbf{w}$ v.s. noise $\boldsymbol\zeta$ & (d) convergence of Alg. \ref{alg:nonlin}
  \end{tabular}
  \caption{Numerical results for Example 4 with $r=0.5$ noise.  In (c), the solid and dashed lines
  refer to the noise $\boldsymbol\zeta$ and the weight $\mathbf{w}$, respectively.}\label{fig:exam4sol}
\end{figure}

The numerical results for Example 4 with data of various noise levels are shown in Table
\ref{tab:exam4}. The accuracy $e$ only deteriorates very mildly as the corruption
percentage $r$ increases from $0.1$ to $0.9$. The result for a typical realization of
noisy data with $r=0.5$ is shown in Fig. \ref{fig:exam4sol}. The solution is not as
accurate as before, since it oscillates slightly around the discontinuities of the true
solution. This is a consequence of the smoothness prior adopted here, which in principle
is unsuited to reconstructing discontinuous profiles. Nonetheless, the solution is
reasonable as the overall profile of the true solution is largely retrieved, and the
magnitude is accurate. The convergence of Algorithm \ref{alg:nonlin} remains very stable
for the discontinuous solution. However, there are several large plateaus, which might be
pruned out by increasing the tolerance $tol$ so as to effect the desired computational
speedup.

\section{Concluding remarks}

In this paper we have developed a robust Bayesian approach to inverse problems subject to
impulsive noises. It explicitly adopts a heavy-tailed $t$ model to cope with data
outliers, and it admits a scale mixture representation, which enables deriving efficient
variational algorithms. The approach has been illustrated on several benchmark linear and
nonlinear inverse problems arising in heat transfer. The numerical results are accurate
and stable even in the presence of a fairly large amount of data outliers, and it is much
more robust compared with the conventional Gaussian model.

There are several avenues deserving further research. We have restricted our attention to
the simplest Markov random field. A natural research problem would be the extension to
general random fields, especially sparsity-promoting prior. Second, it is useful to
develop alternative techniques, e.g.,  based on the Laplace model, and to compare their
merits. The method can only achieve reasonable computational efficiency for medium-scale
problems due to the variance component. It is of interest to develop scalable algorithms
by imposing further restrictions on the approximation. Lastly, rigorous justification of
the excellent performance of the model as well as the algorithm, e.g., consistency and
convergence rate, is of immense theoretical importance, and is to be established.

\section*{Acknowledgements}
This work is supported by Award No. KUS-C1-016-04, made by King Abdullah University of
Science and Technology (KAUST). The author is grateful to two anonymous referees for
their constructive comments, which have led to an improved presentation of the
manuscript.
\bibliographystyle{abbrv}
\bibliography{bayes}

\begin{thebibliography}{10}

\bibitem{AchimBezerianosTsakalides:2001}
A.~Achim, A.~Bezerianos, and P.~Tsakalides.
\newblock Novel {B}ayesian multiscale method for speckle removal in medical
  ultrasound images.
\newblock {\em IEEE Trans. Med. Imag.}, 20(8):772--783, 2001.

\bibitem{AllineyRuzinsky:1994}
S.~Alliney and S.~Ruzinsky.
\newblock An algorithm for the minimization of mixed l1 and l2 norms with
  application to {B}ayesian estimation.
\newblock {\em IEEE Trans. Signal Process.}, 42(3):618--627, 1994.

\bibitem{Attias:2000}
H.~Attias.
\newblock A variational {B}ayesian framework for graphical models.
\newblock In {\em Advances in Neural Information Processing Systems},
  volume~12, pages 209--215, 2000.

\bibitem{Beal:2003}
M.~J. Beal.
\newblock {\em Variational Algorithms for Approximate Bayesian Inference}.
\newblock PhD thesis, Gatsby Computational Neuroscience Unit, University
  College London, 2003.

\bibitem{Beck:1985}
J.~V. Beck, B.~Blackwell, and C.~R. {St. Clair}.
\newblock {\em Inverse {H}eat {C}onduction: {I}ll-{P}osed {P}roblems}.
\newblock Wiley, New York, 1985.

\bibitem{ChanHoNikolova:2005}
R.~H. Chan, C.-H. Ho, and M.~Nikolova.
\newblock Salt-and-pepper noise removal by median-type noise detectors and
  detail-preserving regularization.
\newblock {\em IEEE Trans. Imag. Process.}, 14(10):1479--1485, 2005.

\bibitem{ChappellGrovesWhitcherWoolrich:2009}
M.~A. Chappell, A.~R. Groves, B.~Whitcher, and M.~W. Woolrich.
\newblock Variational {B}ayesian inference for a nonlinear forward model.
\newblock {\em IEEE Trans. Signal Process.}, 57(1):223--236, 2009.

\bibitem{ChoiWette:1969}
S.~C. Choi and R.~Wette.
\newblock Maximum likelihood estimation of the parameters of the {G}amma
  distribution and their bias.
\newblock {\em Technometrics}, 11(4):683--690, 1969.

\bibitem{Clason:2009b}
C.~Clason, B.~Jin, and K.~Kunisch.
\newblock A duality-based splitting method for {$\ell^1$-$TV$} image
  restoration with automatic regularization parameter choice.
\newblock {\em SIAM J. Sci. Comput.}, 32(3):1484--1505, 2010.

\bibitem{ClasonJinKunisch:2010sii}
C.~Clason, B.~Jin, and K.~Kunisch.
\newblock A semismooth {N}ewton method for {$L^1$} data fitting with automatic
  choice of regularization parameters and noise calibration.
\newblock {\em SIAM J. Imaging Sci.}, 3(2):199--231, 2010.

\bibitem{ColliGuerri:1985}
P.~Colli-Franzone, L.~Guerri, S.~Tentoni, C.~Viganotti, S.~Baruffi,
  S.~Spaggiari, and B.~Taccardi.
\newblock A mathematical procedure for solving the inverse potential problem of
  electrocardiography: analysis of the time-space accuracy from in vitro
  experimental data.
\newblock {\em Math. Biosci.}, 77(1--2):353--396, 1985.

\bibitem{DempsterLairdRubin:1977}
A.~P. Dempster, N.~M. Laird, and D.~B. Rubin.
\newblock Maximum likelihood from incomplete data via the {EM} algorithm.
\newblock {\em J. Royal Stat. Soc., Ser. B}, 39(1):1--38, 1977.

\bibitem{Emery:2002}
A.~F. Emery.
\newblock Parameter estimation in the presence of uncertain parameters and with
  correlated data errors.
\newblock {\em Int. J. Thermal Sci.}, 41(6):481--491, 2002.

\bibitem{Emery:2009}
A.~F. Emery.
\newblock Estimating deterministic parameters by {B}ayesian inference with
  emphasis on estimating the uncertainty of the parameters.
\newblock {\em Inv. Probl. Sci. Eng.}, 17(2):263--274, 2009.

\bibitem{Gao:2008}
J.~Gao.
\newblock Robust ${L}^1$ principal component analysis and its {B}ayesian
  variational inference.
\newblock {\em Neur. Comput.}, 20(2):555--572, 2008.

\bibitem{GelmanCarlinSternRubin:2004}
A.~Gelman, J.~B. Carlin, H.~S. Stern, and D.~B. Rubin.
\newblock {\em Bayesian {D}ata {A}nalysis}.
\newblock Chapman \& Hall/CRC, Boca Raton, FL, 2nd edition, 2004.

\bibitem{Geweke:1993}
J.~Geweke.
\newblock Bayesian treatment of the independent student- t linear model.
\newblock {\em J. Appl. Econometrics}, 8(S):S19--40, 1993.

\bibitem{HebertLeahy:1989}
T.~Hebert and R.~Leahy.
\newblock A generalized {EM} algorithm for 3-d {B}ayesian reconstruction from
  {P}oisson data using {G}ibbs priors.
\newblock {\em IEEE Trans. Med. Imag.}, 8(2):194--202, 1989.

\bibitem{Inglese:1997}
G.~Inglese.
\newblock An inverse problem in corrosion detection.
\newblock {\em Inverse Problems}, 13(4):977--994, 1997.

\bibitem{Jin:2008}
B.~Jin.
\newblock Fast {B}ayesian approach for parameter estimation.
\newblock {\em Int. J. Numer. Methods Engrg.}, 76(2):230--252, 2008.

\bibitem{JinLu:2011}
B.~Jin and X.~Lu.
\newblock Numerical identification of a robin coefficient in parabolic
  problems.
\newblock {\em Math. Comput.}, page in press., 2011.

\bibitem{JinZou:2009ip}
B.~Jin and J.~Zou.
\newblock Augmented {T}ikhonov regularization.
\newblock {\em Inverse Problems}, 25(2):025001, 25, 2009.

\bibitem{JinZou:2010jcp}
B.~Jin and J.~Zou.
\newblock Hierarchical {B}ayesian inference for ill-posed problems via
  variational method.
\newblock {\em J. Comput. Phys.}, 229(19):7317--7343, 2010.

\bibitem{JinZou:2010ima}
B.~Jin and J.~Zou.
\newblock Numerical estimation of the {R}obin coefficient in a stationary
  diffusion equation.
\newblock {\em IMA J. Numer. Anal.}, 30(3):677--701, 2010.

\bibitem{JordanGhahramaniJaakkolaSaul:1999}
M.~I. Jordan, Z.~Ghahramani, T.~S. Jaakkola, and L.~K. Saul.
\newblock An introduction to variational methods for graphical models.
\newblock {\em Mach. Learn.}, 37:183--233, 1999.

\bibitem{KaipioSomersalo:2005}
J.~Kaipio and E.~Somersalo.
\newblock {\em Statistical and {C}omputational {I}nverse {P}roblems}.
\newblock Springer-Verlag, New York, 2005.

\bibitem{KaipioKolehmainenSomersaloVauhkonen:2000}
J.~P. Kaipio, V.~Kolehmainen, E.~Somersalo, and M.~Vauhkonen.
\newblock Statistical inversion and {M}onte {C}arlo sampling methods in
  electrical impedance tomography.
\newblock {\em Inverse Problems}, 16(5):1487--1522, 2000.

\bibitem{Koutsourelakis:2009}
P.~S. Koutsourelakis.
\newblock A multi-resolution, non-parametric, {B}ayesian framework for
  identification of spatially-varying model parameters.
\newblock {\em J. Comput. Phys.}, 228(17):6184--6211, 2009.

\bibitem{LangeLittleTaylor:1989}
K.~L. Lange, R.~J.~A. Little, and J.~M.~G. Taylor.
\newblock Robust statistical modeling using the {$t$} distribution.
\newblock {\em J. Amer. Statist. Assoc.}, 84(408):881--896, 1989.

\bibitem{MaZabaras:2009}
X.~Ma and N.~Zabaras.
\newblock An efficient {B}ayesian inference approach to inverse problems based
  on an adaptive sparse grid collocation method.
\newblock {\em Inverse Problems}, 25(3):035013, 27, 2009.

\bibitem{MarzoukXiu:2009}
Y.~Marzouk and D.~Xiu.
\newblock A stochastic collocation approach to {B}ayesian inference in inverse
  problems.
\newblock {\em Comm. Comput. Phys.}, 6(4):826--847, 2009.

\bibitem{MarzoukNajmRahn:2007}
Y.~M. Marzouk, H.~N. Najm, and L.~A. Rahn.
\newblock Stochastic spectral methods for efficient bayesian solution of
  inverse problems.
\newblock {\em J. Comput. Phys.}, 224(2):560--586, 2007.

\bibitem{OsmanBeck:1988}
A.~M. Osman and J.~V. Beck.
\newblock Nonlinear inverse problem for the estimation of time-and-space
  dependent heat transfer coefficients.
\newblock {\em J. Thermophys.}, 3:146--152, 1988.

\bibitem{SambridgeMosegaard:2002}
M.~Sambridge and K.~Mosegaard.
\newblock Monte {C}arlo methods in geophysical inverse problems.
\newblock {\em Rev. Geophys.}, 40(3):3--29, 2002.

\bibitem{SatoTaku:2004}
M.~Sato, T.~Yoshioka, S.~Kajihara, K.~Toyama, N.~Goda, K.~Doya, and M.~Kawato.
\newblock Hierarchical {B}ayesian estimation for {MEG} inverse problem.
\newblock {\em NeuroImage}, 23(3):806--826, 2004.

\bibitem{SuHewitt:2004}
J.~Su and G.~F. Hewitt.
\newblock Inverse heat conduction problem of estimating time-varying heat
  transfer coefficient.
\newblock {\em Numer. Heat Transfer, Part A}, 45(8):777--789, 2004.

\bibitem{Tarantola:2005}
A.~Tarantola.
\newblock {\em Inverse {P}roblem {T}heory and {M}ethods for {M}odel {P}arameter
  {E}stimation}.
\newblock SIAM, Philadelphia, PA, 2005.

\bibitem{Tipping:2005}
M.~E. Tipping and N.~D. Lawrence.
\newblock Variational inference for student-t models: robust bayesian
  interpolation and generalised component analysis.
\newblock {\em Neurocomput.}, 69(1--3):123--141, 2005.

\bibitem{WangZabaras:2004ijhmt}
J.~Wang and N.~Zabaras.
\newblock A {B}ayesian inference approach to the inverse heat conduction
  problem.
\newblock {\em Int. J. Heat Mass Transfer}, 47(17-18):3927--3941, 2004.

\bibitem{WangZabaras:2005ip}
J.~Wang and N.~Zabaras.
\newblock Hierarchical {B}ayesian models for inverse problems in heat
  conduction.
\newblock {\em Inverse Problems}, 21(1):183--206, 2005.

\end{thebibliography}

\end{document}